\documentclass[11pt]{amsart}
\usepackage{amsmath,upgreek}
\usepackage[arrow,matrix,curve]{xy}

\theoremstyle{plain}
\newtheorem{thm}{Theorem}[section]
\newtheorem{lm}[thm]{Lemma}
\newtheorem{cor}[thm]{Corollary}
\newtheorem{pro}[thm]{Proposition}
\theoremstyle{definition}
\newtheorem{df}[thm]{Definition}

\def\sus{\mathsf s}
\def\bb{\mathsf b}
\def\dd{\mathsf d}
\def\mm{\mathsf m}
\def\ff{\mathsf f}
\def\gg{\mathsf g}
\def\hh{\mathsf h}
\def\pp{\mathsf p}
\def\qq{\mathsf q}
\def\FF{\mathsf F}
\def\GG{\mathsf G}
\def\HH{\mathsf H}
\def\QQ{\mathsf Q}
\def\TT{\uptau}

\newcommand\al{\alpha}
\newcommand\gam{\gamma}
\newcommand\sig{\sigma}
\newcommand\veps{\varepsilon}
\newcommand\Del{\mathit{\Delta}}
\def\Al{{\boldsymbol{\alpha}}}

\newcommand\aaa{{\boldsymbol{a}}}
\newcommand\bbb{{\boldsymbol{b}}}
\newcommand\ccc{{\boldsymbol{c}}}
\newcommand\hhh{{\boldsymbol{h}}}
\newcommand\xxx{{\boldsymbol{x}}}
\newcommand\yyy{{\boldsymbol{y}}}
\newcommand\ooo{{\boldsymbol{0}}}
\newcommand\kk{{\boldsymbol{k}}}
\newcommand\Fi{{\mathrm{F}}}
\newcommand\Gr{{\mathrm{Gr}}}
\newcommand\Coh{{\mathrm{H}}}
\newcommand\Z{\mathbb{Z}}

\def\op{\oplus}
\def\ox{\otimes}
\def\bd{\partial}
\def\eq{&\!\!=\!\!&}
\def\vc{&\!\!&}
\def\la{\langle}
\def\ra{\rangle}
\def\bart{\overline{T}}

\newcommand\one{\mathsf 1}
\def\shuffle{\Join}
\newcommand\mat[1]{{\mathrm{Mat}}^{(#1)}}
\newcommand\pqr{_r^{pq}}
\newcommand\opqr{_{r-1}^{p+1,q-1}}
\newcommand\pqrr{_r^{p+r,q-r+1}}
\newcommand\set[2]{\{#1\,|\,#2\}}
\newcommand\two[4]{\left\{    \begin{array}{ll}
                                        {#1}, & {\mbox{if }} {#2}, \\
                                        {#3}, & {\mbox{if }} {#4}
                                        \end{array}     \right.}
\newcommand\ubrace[2]{\underbrace{{#1},\dots,{#1}}_{#2\;\text{times}}}
\newcommand\dbrace[2]{\underbrace{{#1}\odot\cdots\odot{#1}}_{#2\;\text{times}}}
\newcommand\pbrace[2]{\underbrace{{#1}\op\dots\op{#1}}_{#2\;\text{times}}}
\newcommand\xbrace[2]{\underbrace{{#1}\ox\cdots\ox{#1}}_{#2\;\text{times}}}

\begin{document}

$\,$

\baselineskip.525cm

\vspace{-1.65cm}

\begin{flushright}
{\tt arXiv:0912.1775v2[math.AT]}\\
December, 2009 (revised)
\end{flushright}

\vspace{1cm}

\title[A-infinity algebras]
{Massey Product and Twisted Cohomology\\ of $\boldsymbol{A}$-infinity Algebras}
\author[Weiping Li]{Weiping Li}
\address{Department of Mathematics, Oklahoma State University, 
Stillwater, OK 74078-0613, USA}
\email{wli@math.okstate.edu}
\author[Siye Wu]{Siye Wu}
\address{Department of Mathematics, University of Colorado, 
Boulder, CO 80309-0395, USA and \newline
\hspace*{.175in}Department of Mathematics,
The University of Hong Kong, Pokfulam, Hong Kong}
\email{swu@math.colorado.edu, swu@maths.hku.hk}
\date{\today}


\begin{abstract}
We study the twisted cohomology groups of $A_\infty$-algebras defined by 
twisting elements and their behavior under morphisms and homotopies using
the bar construction.
We define higher Massey products on the cohomology groups of general 
$A_\infty$-algebras and establish the naturality under morphisms and their
dependency on defining systems.
The above constructions are also considered for $C_\infty$-algebras.
We construct a spectral sequence converging to the twisted cohomology groups
and show that the higher differentials are given by the $A_\infty$-algebraic
Massey products.

\end{abstract}
\maketitle

\section{Introduction}\label{intro}

The concept of $A_\infty$-algebras was introduced by Stasheff \cite{St1,St2}
for studying multiplication operations which satisfy associativity up to
homotopy. 
Since then it has played a crucial role in homotopy theory.
The $A_\infty$-structure on the cohomology of a topological space determines
the cohomology of its loop space and can be applied to the cohomology of
fiber bundles \cite{K80}. 
Moreover, it determines the rational homotopy type of $1$-connected spaces
\cite{K08}. 
Recently, the subject finds applications in many areas of algebra, topology,
geometry and mathematical physics, including homological mirror symmetry
\cite{KS}.

Motivated by the work on twisted cohomology of the de Rham complex
\cite{RW,AS,MW,LLW}, we study the twisted cohomology of $A_\infty$-algebras.
In addition, we define higher Massey products on the cohomology of a 
general $A_\infty$-algebra with a possibly non-associative multiplication.
We then construct a spectral sequence converging to the twisted cohomology
and relate the higher differentials to the higher Massey products.

The paper is organized as follows.
In Section~\ref{rev}, we recall the notions of $A_\infty$-algebras,
morphisms of $A_\infty$-algebras and homotopies of morphisms.
We express these concepts using the bar construction.
The $A_\infty$-structure on the cohomology group is also described.
The $C_\infty$-algebras are discussed as a special case. 
In Section~\ref{twist}, we study the twisted cohomology group of a
differential on the $A_\infty$-algebra deformed by a twisting element.
(These concepts simplify when the $A_\infty$-algebra is $C_\infty$.)
We put special emphasis on the use of the bar construction, which clarifies
many concepts and calculations. 
The twisted cohomology groups are naturally isomorphic if the twisting
elements are homotopic equivalent.
We show that a morphism of $A_\infty$-algebras maps a twisting element to
another and induces a homomorphism on the twisted cohomology groups. 
We also find the relation of the induced homomorphisms from two homotopic
morphisms of $A_\infty$-algebras.
In Section~\ref{massey}, we introduce the triple and higher Massey products
on the cohomology of a general $A_\infty$-algebra.
A crucial ingredient is the matric $A_\infty$-algebra formed by matrices of
elements in the original $A_\infty$-algebra, together with its properties.
We introduce an equivalence relation on the set of defining systems under
which the (higher) Massey products take the same value; this is of
interest even in the classical case.
We establish some properties of the Massey product and the naturality under
morphisms of $A_\infty$-algebras.
In particular, we clarify and generalize the folklore relation of the Massey
product and the $A_\infty$-structure on the cohomology (of differential graded
algebras) to the context of $A_\infty$-algebras.
We also study the Massey products for $C_\infty$-algebras.
In Section~\ref{spectral}, we construct a spectral sequence associated to a
natural filtration on a $\Z$-graded $A_\infty$-algebra (assuming the twisting
element is of positive degree) that converges to the $\Z_2$-graded twisted
cohomology group.
We give a complete description of the higher differentials in terms of the
higher Massey products of the $A_\infty$-algebra.
We show that a morphism of $A_\infty$-algebras induces a morphism of spectral
sequences.
This result is applied to the quasi-isomorphism of the cohomology group and
the original $A_\infty$-algebra.
In the Appendix, we present a construction of spectral sequence that is
slightly different from what we can find in the literature but suits the
purpose for the previous section.

\section{$A_\infty$-algebras and morphisms}\label{rev}

In this section, we recall the definitions of $A_\infty$- and 
$C_\infty$-algebras, their morphisms, homotopy of morphisms,
and the bar construction.

\subsection{$A_\infty$-algebras, morphisms and homotopy}

\begin{df}\label{df-b}
Let $\kk$ be a field.
An {\em $A_{\infty}$-algebra} $(A,\{\bb_n\})$ over $\kk$ is a $\Z$- or 
$\Z_2$-graded vector space $A=\bigoplus_pA^p$ over $\kk$ with graded
homogeneous $\kk$-multilinear maps $\bb_n\colon A^{\ox n}\to A$ ($n\geq1$)
of degree $1$ satisfying 
\[ \sum_{\substack{r,t\ge0,\,s\ge1\\r+s+t=n}}
          \bb_{r+1+t}\circ(\one^{\ox r}\ox\bb_s\ox\one^{\ox t})=0,\]
where $\one=\one_A$ is the identity map on $A$.
\end{df}

When maps are evaluated on elements, we follow Koszul's sign rule 
$(\ff\ox\gg)(x\ox y)=(-1)^{|\gg||x|}\ff(x)\ox\gg(y)$, where $\ff$, $\gg$ 
are graded homogeneous maps of degrees $|\ff|$, $|\gg|$ and $x$, $y$ are
homogeneous elements of degrees $|x|$, $|y|$, respectively.
So the above identity on $a_1,\dots,a_n\in A$ is
\[  \sum_{\substack{r,t\ge0,\,s\ge1\\r+s+t=n}}
            \bb_{r+1+t}(\bar a_1,\dots,\bar a_r,\bb_s(a_{r+1},\dots,a_{r+s}),
                              a_{r+s+1},\dots,a_n)=0, \]
where $\bar a=(-1)^{|a|}\,a$ for any homogeneous element $a\in A$.
For $n=1$, the identity is $\bb_1\circ\bb_1=0$, i.e., $(A^\bullet,\bd)$,
where $\bd=\bb_1$, is a cochain complex.
For $n\ge2$, the meaning of the identities is best seen from the desuspended
version. 
The $A_\infty$-algebras defined here differ from the original ones (see
\cite{Pr,Ke} for reviews) by a suspension.
Let $\sus\colon\sus^{-1}A\to A$ be the suspension map of degree $-1$ given by
the identification $(\sus^{-1}A)^p=A^{p-1}$.
Then the multilinear maps $\mm_n=\pm\,\sus^{-1}\circ\bb_n\circ\sus^{\ox n}
\colon(\sus^{-1}A)^{\ox n}\to\sus^{-1}A$ (each with an appropriate sign)
satisfy a similar set of identities.
The identity for $n=2$ is the Leibniz property for the product $\mm_2$ while
that for $n=3$ says that $\mm_2$ is associative up to a homotopy given by
$\mm_3$.
If $\mm_n=0$ (or equivalently, $\bb_n=0$) for all $n\ge3$, then $A$ is 
a differential graded algebra.

\begin{df}\label{df-f}
Let $(A,\{\bb^A_n\})$, $(B,\{\bb^B_n\})$ be $A_\infty$-algebras.
A {\em morphism} $\ff\colon A\to B$ of $A_\infty$-algebras is a family 
$\ff=\{\ff_n\}_{n\ge1}$ of graded homogeneous $\kk$-multilinear maps
$\ff_n\colon A^{\ox n}\to B$ of degree $0$ satisfying, for each $n\ge1$,
the identity
\[\sum_{\substack{r,t\ge0,\,s\ge1\\r+s+t=n}}
         \ff_{r+1+t}\circ(\one^{\ox r}\ox\bb^A_s\ox\one^{\ox t})
=\sum_{\substack{r\ge0;\,i_1,\dots,i_r>0\\i_1+\cdots+i_r=n}}
         \bb^B_r(\ff_{i_1}\ox\cdots\ox\ff_{i_r}).                           \]
\end{df}

For $n=1$, the above identity means that $\ff_1\colon(A,\bd^A)\to(B,\bd^B)$
is a morphism of cochain complexes.
For $n=2$, $\ff_1$ commutes with the operations $\mm^A_2,\mm^B_2$ (in the
desuspended version) up to a homotopy given by $\ff_2$.
A morphism $\ff$ is called a {\em quasi-isomorphism} if $\ff_1$ is a 
quasi-isomorphism.
It is {\em strict} if $\ff_n=0$ for all $n\ge2$.
The composition of two $A_\infty$-algebra morphisms $\ff\colon A\to B$ and
$\gg\colon B\to C$ is given by, for any $n\ge1$,
\[ (\gg\circ\ff)_n=\sum_{\substack{r\ge0;\,i_1,\dots,i_r>0\\i_1+\cdots+i_r=n}}
                         \gg_r\circ(\ff_{i_1}\ox\cdots\ox\ff_{i_r}).      \]
In particular, $(\gg\circ\ff)_1=\gg_1\circ\ff_1$.

\begin{df}\label{df-h}
Let $(A,\{\bb^A_n\})$, $(B,\{\bb^B_n\})$ be $A_\infty$-algebras.
Two morphisms $\{\ff_n\}$, $\{\gg_n\}$ from $A$ to $B$ are {\em homotopic}
(through $\{\hh_n\}$) if there exists a family of graded homogeneous 
$\kk$-multilinear maps $\hh_n\colon A^{\ox n}\to B$ ($n\ge1$) of degree $-1$
such that for any $n\ge 1$,
\[ \gg_n-\ff_n=\sum\bb^B_{r+1+t}\circ
   (\gg_{i_1}\ox\cdots\ox\gg_{i_r}\ox\hh_s\ox\ff_{j_1}\ox\cdots\ff_{j_t})
   +\sum\hh_{r+1+t}\circ(\one^{\ox r}\ox\bb^A_s\ox\one^{\ox t}),    \]
where the first sum runs over $r,t\ge0$, $s\ge1$, 
$i_1+\cdots+i_r+j_1+\cdots+j_t=n$ ($i_1,\dots,i_r,j_1,\dots,j_t>0$) and
the second runs over $r,t\ge0$, $s\ge1$, $r+s+t=n$.
\end{df}

Homotopy is an equivalence relation on the set of morphisms.\footnote{See for
example \cite{Pr,Hi,LH}. We thank B.\ Keller for providing the references.} 

\subsection{The bar construction}

A more conceptual way of describing $A_\infty$-structures is through the use
the bar construction \cite{St2}.
Recall that the reduced tensor coalgebra on a vector space $A$ over $\kk$ is
$\bart A=\bigoplus_{n\ge1} A^{\ox n}$.
There is a comultiplication
$\Del\colon\bart A\to\bart A\ox\bart A$ given by
\[  \Del\,[a_1\ox\cdots\ox a_n]
=\sum_{r=1}^{n-1}\,[a_1\ox\cdots\ox a_r]\ox[a_{r+1}\ox\cdots\ox a_n],   \]
where $0<r<n$ and $a_1,\dots,a_n\in A$.
For an $A_\infty$-algebra $A$, the map $\sum_{n\ge1}\bb_n\colon\bart A\to A$
lifts uniquely to a graded coderivation $\bb\colon\bart A\to\bart A$ 
of degree $1$ satisfying
\[ \Del\circ\bb=(\one\ox\bb+\bb\ox\one)\circ\Del. \]
The conditions that $(A,\{\bb_n\})$ is an $A_\infty$-algebra is equivalent
to the identity $\bb\circ\bb=0$, i.e., $\bb$ is a coalgebra differential
on $\bart A$ of degree $1$.

Similarly, a collection of $\kk$-multilinear maps
$\ff=\{\ff_n\colon A^{\ox n}\to B\}_{n\ge1}$ defines a map
$\sum_{n\ge1}\ff_n\colon\bart A\to B$ which lifts uniquely to a coalgebra
morphism $\FF\colon\bart A\to\bart B$ of degree $0$.
If $(A,\{\bb^A_n\})$, $(B,\{\bb^B_n\})$ are $A_\infty$-algebras, the
condition that $\ff=\{\ff_n\}\colon A\to B$ is an $A_\infty$-algebra
morphism is equivalent to $\FF\circ\bb^A=\bb^B\circ\FF$, i.e.,
$\FF\colon(\bart A,\bb^A)\to(\bart B,\bb^B)$ is a morphism of graded
differential coalgebras.
If $\gg\colon B\to C$ is another morphism of $A_\infty$-algebras, the
composition of morphisms $\gg\circ\ff\colon A\to C$ corresponds to the
usual composition $\GG\circ\FF\colon\bart A\to\bart C$.

Two $A_\infty$-algebra morphisms $\ff,\gg\colon A\to B$ are homotopic if and
only if $\FF,\GG\colon\bart A\to\bart B$ are homotopic as morphisms of graded
differential coalgebras, i.e., there is a map $\HH\colon\bart A\to\bart B$
of degree $-1$ such that
\[ \GG-\FF=\bb^B\circ\HH+\HH\circ\bb^A,\quad 
\Del^B\circ\HH=(\GG\ox\HH+\HH\ox\FF)\circ\Del^A.    \]

\subsection{$C_\infty$-algebras and morphisms}\label{C-alg-mor}

If $A$ is a vector space over $\kk$, the shuffle product \cite{Mac} on
$\bart A$ is given by, for any $0<r<n$ and $a_1\dots,a_n\in A$,
\[  [a_1\ox\cdots\ox a_r]\shuffle[a_{r+1}\ox\cdots\ox a_n]=\sum_{\sig\in
        S_{r,n}}(-1)^{\veps(\sig)}\,a_{\sig(1)}\ox\cdots\ox a_{\sig(n)},\]
where $S_{r,n}$ is the set of permutations $\sig\in S_n$ on the set
$\{1,\dots,n\}$ such that $\sig(i)<\sig(j)$ if either $1\le i<j\le r$ or
$r+1\le i<j\le n$ and $\veps(\sig)=\sum |a_i||a_j|$, summing over the pairs
$(i,j)$ with $1\le i\le r<j\le n$ and $\sig(i)>\sig(j)$.
This product is graded commutative and associative on $\bart A$, making it,
together with the comultiplication $\Del$, a graded bi-algebra.

\begin{df}\label{cinf}
A {\em $C_\infty$-algebra} $(A,\{\bb_n\})$ is an $A_\infty$-algebra such that
for each $n\ge2$, $\bb_n=0$ on $\bart A\shuffle\bart A$.
A {\em morphism} $\ff=\{\ff_n\}_{n\ge1}\colon A\to B$ of $C_\infty$-algebras
is a morphism of $A_\infty$-algebras such that for each $n\ge2$, $\ff_n=0$ on
$\bart A\shuffle\bart A$.
Two $C_\infty$-algebra morphisms $\ff,\gg\colon A\to B$ are {\em homotopic}
if they are homotopic as $A_\infty$-algebra morphisms through $\{\hh_n\}$
such that for each $n\ge2$, $\hh_n=0$ on $\bart A\shuffle\bart A$. 
\end{df}

For $n=2$, the condition $\bb_2(A\shuffle A)=0$ is equivalent to the graded
commutativity of the product $\mm_2$ on $\sus^{-1}A$.
Thus, $C_\infty$-algebras are ``graded commutative'' $A_\infty$-algebras
(see \cite{K08} for a history of and references on $C_\infty$-algebras).
The conditions on $C_\infty$-structures and morphisms can be stated concisely
using the bar construction \cite{K08}.
The $A_\infty$-structure on $A$ is a $C_\infty$-structure if and only if $\bb$
is a derivation on the graded bi-algebra $\bart A$, i.e., for any
$x,y\in\bart A$,
\[   \bb(x\shuffle y)=\bb(x)\shuffle y+\bar x\shuffle\bb(y).  \]
Similarly, an $A_\infty$-morphism $\ff\colon A\to B$ is $C_\infty$ if and only
if $\FF\colon\bart A\to\bart B$ if a morphism of bi-algebras, i.e., for any
$x,y\in\bart A$,
\[   \FF(x\shuffle y)=\FF(x)\shuffle\FF(y).   \]
Using the same method, one can further show that if $\hh=\{\hh_n\}$
is a homotopy of $C_\infty$-algebra morphisms $\ff,\gg\colon A\to B$ if and
only if for any $x,y\in\bart A$,
\[ \HH(x\shuffle y)=\GG(\bar x)\shuffle\HH(y)+\HH(x)\shuffle\FF(y). \] 

\subsection{$A_\infty$-structure on the cohomology group}\label{sec-ha}
Given an $A_\infty$-algebra $(A,\{\bb_n\})$, the cohomology 
$\Coh(A)=\Coh(A,\bd)$ is an associative graded algebra
under the operation $\bar\bb_2$ induced by $\bb_2$.
In fact, the cohomology $\Coh(A)$ has an $A_\infty$-algebra structure
$\{\bar\bb_n\}$ with $\bar \bb_1=0$ \cite{K80,K82,Me,KS}.
There is a quasi-isomorphism of $A_\infty$-algebras $\Coh(A)\to A$ lifting 
the identity map of $\Coh(A)$.
Such an $A_\infty$-structure on $\Coh(A)$ is unique up to isomorphisms of
$A_\infty$-algebras.

We describe the $A_\infty$-algebra $(\Coh(A),\{\bar\bb_n\})$ and
the quasi-isomorphism $\qq\colon \Coh(A)\to A$.
Let $\pp_1\colon A\to \Coh(A)$ be a $\kk$-linear map which sends any
closed element to the cohomology class it represents; $\pp_1$ is determined
by the choice of a subspace in $A$ that is transverse to $\ker\bd$.
Let $\qq_1\colon \Coh(A)\to\ker\bd\subset A$ be a $\kk$-linear map
such that $\pp_1\circ\qq_1=\one_{\Coh(A)}$.
Then there is a homogeneous $\kk$-linear map $\hh_1\colon A\to A$ of degree
$-1$ such that $\one_A-\qq_1\circ\pp_1=\bd\circ\hh_1+\hh_1\circ\bd$.
The $A_\infty$-structure $\{\bar\bb_n\}$ on $\Coh(A)$ and the
quasi-isomorphism $\qq=\{\qq_n\}\colon\Coh(A)\to A$ can be expressed
explicitly as a sum over the oriented rooted planar trees \cite{KS}.
Alternatively, they can be obtained inductively by \cite{CG}
\begin{eqnarray*}
\bar\bb_n\eq\sum_{r=2}^n\sum_{\substack{i_1+\cdots+i_r=n\\
  i_1,\dots,i_r>0}}\pp_1\circ\bb_r\circ(\qq_{i_1}\ox\cdots\ox\qq_{i_r}),  \\
\qq_n\eq\sum_{r=2}^n\sum_{\substack{i_1+\cdots+i_r=n\\ i_1,\dots,i_r>0}}
          \hh_1\circ\bb_r\circ(\qq_{i_1}\ox\cdots\ox\qq_{i_r}),
\end{eqnarray*}
respectively, for any $n\ge2$.
When $A$ is a differential graded algebra, the inductive formulas simplify
to \cite{LPWZ}
\[ \bar\bb_n=\sum_{i=1}^{n-1}\pp_1\circ\bb_2\circ(\qq_i\ox\qq_{n-i}),\quad
\qq_n=\sum_{i=1}^{n-1}\hh_1\circ\bb_2\circ(\qq_i\ox\qq_{n-i}).    \]
Finally, if $(A,\{\bb_n\})$ is a $C_\infty$-algebra, then so is
$(\Coh(A),\{\bar\bb_n\})$ \cite{CG}.

\section{Twisting elements and twisted cohomology}\label{twist}

\subsection{Twisting elements and twisted cohomology of $A_\infty$-algebras}
\label{sec-tw}

We construct a deformed differential on the $A_\infty$-algebra and study its
cohomology group.
Let $(A,\{\bb_n\})$ be an $A_\infty$-algebra.
We fix an element $h\in A$.
Define a map $\TT_h\colon A\to\bart A$ by
\[ a\in A\mapsto\TT_h(a)=\sum_{n=0}^\infty\xbrace{h}{n}\ox\,a
           =a+h\ox a+h\ox h\ox a+\cdots. \]
$\TT_h$ is a $\kk$-linear isomorphism from $A$ onto its image, and the
inverse map is given by $\TT_h^{-1}(x)=x-h\ox x$ for $x\in\TT_h(A)$.
We set $\TT_h(1)=1+\TT_h(h)\in\kk\op\bart A$.
Then $\TT_h(a)=\TT_h(1)\ox a$.

It is clear that $\Del(\TT_h(a))=\TT_h(h)\ox\TT_h(a)$ for any $h,a\in A$.
In particular, $\Del(\TT_h(h))=\TT_h(h)\ox\TT_h(h)$, i.e., $\TT_h(h)$ can
be regarded as a morphism of coalgebras from $\kk$ (with the obvious
comultiplication) to $\bart A$.

We define a twisted differential $\bd_h\colon A\to A$ by\footnote{During
the preparation of the paper, we learned that the differential was previously
introduced and studied in \cite{Pr}.
Our approach below relies instead on the bar construction. In the case when
$A$ is a $C_\infty$-algebra, which is treated in \S\ref{C-twist}, the same
differential was also obtained recently by E.\ Getzler.}
\[ a\in A\mapsto\bd_ha=\sum_{n=0}^\infty\bb_{n+1}(\ubrace{h}{n},a)=
   \bd a+\bb_2(h,a)+\bb_3(h,h,a)+\cdots.  \]

\begin{lm}\label{bdh}
If $h\in A$, then
\[ \bb\circ\TT_h=\TT_{\bar h}\circ\bd_h+\TT_{\bar h}(\bd_hh)\ox\TT_h, \]
where $\bar h=(-1)^{|h|}h$.
\end{lm}
\noindent{\em Proof:}
By the definition of $\bb$, we have, for any $a\in A$,
\begin{eqnarray*}
\bb(\TT_h(a))
\eq\sum_{n=0}^\infty\bb(\xbrace{h}{n}\ox\,a)                \\
\eq\sum_{n=0}^\infty\left(
 \sum_{\substack{r,s\ge0\\r+s=n}}\xbrace{\bar h}{r}\ox\,
       \bb_{s+1}(\ubrace{h}{s},a)\;\;+\right.              \\
& &\left.\quad+\sum_{\substack{r,t\ge0,s\ge1\\r+s+t=n}}
  \xbrace{\bar h}{r}\ox\,\bb_s(\ubrace{h}{s})\ox\xbrace{h}{t}\ox\,a\right) \\
\eq\TT_{\bar h}(1)\ox\bd_ha+\TT_{\bar h}(1)\ox\bd_hh\ox\TT_h(a)   \\
\eq\TT_{\bar h}(\bd_ha)+\TT_{\bar h}(\bd_hh)\ox\TT_h(a).
\end{eqnarray*}
\qed

\begin{cor}\label{twisting}
(i) For any $h\in A$,
\[ \sum_{r,t\ge0}\bb_{r+1+t}(\ubrace{\bar h}{r},\bd_hh,\ubrace{h}{t})=0; \]

(ii) $\bd_hh=0$ if and only if $\bb(\TT_h(h))=0$;

(ii) For any $h,a\in A$,
\[ \bd_{\bar h}\bd_ha=
   -\sum_{r,t\ge0}\bb_{r+t+2}(\ubrace{\bar h}{r},\bd_hh,\ubrace{h}{t},a).  \]

\end{cor}
\noindent{\em Proof:}
Taking $a=h$ in the proof of Lemma~\ref{bdh}, we get 
\[     \bb(\TT_h(h))=\TT_{\bar h}(1)\ox\bd_hh\ox\TT_h(1),   \]
from which (ii) follows.
Applying $\bb$ to both sides and projecting to $A$, we get (i).
Applying $\bb$ to the identity in Lemma~\ref{bdh}, we get
\[ \TT_h(\bd_{\bar h}\bd_ha)=-\TT_h(\bd_{\bar h}\bar h)\ox\TT_{\bar h}(\bd_ha)
                             -\bb(\TT_{\bar h}(\bd_hh)\ox\TT_h(a)). \]
(iii) follows by a projection onto $A$.
\qed

\begin{df}
If $A$ is an $A_\infty$-algebra, an element $h\in A$ is a
{\em twisting element} if it is of even degree and $\bd_hh=0$.
\end{df}

When $A$ is a differential graded algebra, the condition $\bd_hh=0$ reduces
to the Maurer-Cartan equation $\bd h+\bb_2(h,h)=0$.
In general, $\bb(\TT_h(h))=0$ means that $\TT_h(h)$ is a morphism of
differential graded coalgebras from $\kk$ (with the trivial coderivation)
to $\bart A$.
So $\TT_h(h)\in\mathrm{Hom}(\kk,\bart A)$ is a twisted cochain in the sense
of Brown \cite{Br,Pr}.
We refer the readers to \cite{St09} for a history on twisting elements and
twisting cochains, and to the references therein as well as \cite{Pr,K86}.

\begin{thm}\label{BDh} 
If $h$ is a twisting element of an $A_\infty$-algebra $A$, then

(i) $\bb\circ\TT_h=\TT_h\circ\bd_h$ on $A$;

(ii) $\bb$ preserves the subspace $\TT_h(A)\subset\bart A$;

(iii) $\bd_h^2=0$ on $A$.
\end{thm}
\noindent{\em Proof:}
(i) follows from Lemma~\ref{bdh} since $\bar h=h$ and $\bd_hh=0$;\\
(ii) follows easily from (i);\\
(iii) follows from $\bd_h=\TT_h^{-1}\circ\bb\circ\TT_h$ and $\bb^2=0$.
\qed

\begin{df}
If $h$ is a twisting element of the $A_\infty$-algebra $A$, the {\em twisted
cohomology} of $A$ (twisted by $h$) is $\Coh_h(A)=\Coh(A,\bd_h)$.
\end{df}

We note that $\sus^{-1}h$ is of odd degree if $h$ is of even degree.
Although $A$ can be either $\Z$- or $\Z_2$-graded, the twisted cohomology
$\Coh(A)$ is always $\Z_2$-graded.
A special case is the cohomology of the de Rham complex twisted by
a closed form of odd degree \cite{RW,AS,MW,LLW}.

\subsection{Induced morphisms on twisted cohomology groups}

Suppose $\ff=\{\ff_n\}\colon A\to B$ is a morphism of $A_\infty$-algebras,
inducing $\FF\colon\bart A\to\bart B$ on the bar construction.
For any $h\in A$, define $\FF_h\colon A\to B$ by
\[ a\in A\mapsto\FF_h(a)=\sum_{n=0}^\infty\ff_{n+1}(\ubrace{h}{n},a).  \]

\begin{lm}\label{morh}
For any $h\in A$, we have the identity
\[ \FF\circ\TT^A_h=\TT^B_{\FF_h(h)}\circ\FF_h. \]
In particular,
$\FF\colon\TT^A_h(A)\subset\bart A\to\TT^B_{\FF_h(h)}(B)\subset\bart B$
and $\FF(\TT^A_h(h))=\TT^B_{\FF_h(h)}(\FF_h(h))$.
\end{lm}

\noindent{\em Proof:}
For any $a\in A$, by the definition of $\FF$, we have
\begin{eqnarray*}
\FF(\TT^A_h(a))
\eq\sum_{n=0}^\infty\FF(\xbrace{h}{n}\ox\,a)  \\
\eq\sum_{n=0}^\infty\sum_{r=0}^n
   \sum_{\substack{i_0\ge0,\;i_1\dots,i_r\ge1\\i_0+\cdots+i_r=n}}
   \ff_{i_r}(\ubrace{h}{i_r})\ox\cdots\ox\ff_{i_1}(\ubrace{h}{i_1})
   \ox\ff_{i_0+1}(\ubrace{h}{i_0},a)                         \\
\eq\sum_{r=0}^\infty\xbrace{\FF_h(h)}{r}\ox\,\FF_h(a)        \\
\eq\TT^B_{\FF_h(h)}(\FF_h(a)).
\end{eqnarray*}
The rest follows easily.
\qed

\begin{thm}\label{MORh}
Suppose $\ff\colon A\to B$ is a morphism of $A_\infty$-algebras and $h$ is
a twisting element in $A$.
Then

(i) $\FF_h(h)$ is a twisting element in $B$;

(ii) $\FF_h\colon(A,\bd_h)\to(B,\bd_{\FF_h(h)})$ is a cochain map, i.e.,
$\FF_h\circ\bd^A_h=\bd^B_{\FF_h(h)}\circ\FF_h$;

(iii) if $\gg\colon B\to C$ is another morphism of $A_\infty$-algebras, then
$(\GG\circ\FF)_h=\GG_{\FF_h(h)}\circ\FF_h$.
\end{thm}
\noindent{\em Proof:}
By Theorem~\ref{BDh}(i), Lemma~\ref{morh} and by 
$\FF\circ\bb^A=\bb^B\circ\FF$, we have
\begin{eqnarray*}
\TT^B_{\FF_h(h)}\circ\bd^B_{\FF_h(h)}\circ\FF_h
\eq\bb^B\circ\TT^B_{\FF_h(h)}\circ\FF_h=F\circ\bb^A\circ\TT^A_h
=\bb^B\circ\FF\circ\TT^A_h                                            \\
\eq\FF\circ\bb^A\circ\TT^A_h
=\FF\circ\TT^A_h\circ\bd^A_h=\TT^B_{\FF_h(h)}\circ\FF_h\circ\bd^A_h.
\end{eqnarray*}
(ii) follows from the injectivity of $\TT^B_{\FF_h(h)}$ whereas (i) is
from $\bd^B_{\FF_h(h)}\FF_h(h)=\FF_h(\bd^A_hh)=0$.\\
Next, we observe that, for any $a\in A$,
\[ \GG(\FF(\TT^A_h(a)))=\GG(\TT^B_{\FF_h(h)}(\FF_h(a)))=
   \TT^C_{\GG_{\FF_h(h)}(\FF_h(h))}(\GG_{\FF_h(h)}(\FF_h(a))).    \]
On the other hand, we have
\[  \GG(\FF(\TT^A_h(a)))=\TT^C_{(\GG\circ\FF)_h(h)}((\GG\circ\FF)_h(a)).  \]
By setting $a=h$, we obtain $(\GG\circ\FF)_h(h)=\GG_{\FF_h(h)}(\FF_h(h))$
and (iii) follows.
\qed

\begin{cor}
Under the above assumptions,

(i) there is an induced homomorphism
$(\FF_h)_*\colon \Coh_h(A)\to \Coh_{\FF_h(h)}(B)$;

(ii) there is a commutative diagram
\[ \xymatrix@=1pc{& \Coh_{\FF_h(h)}(B)\ar[dr]^{(\GG_{\FF_h(h)})_*} &  \\
   \Coh_h(A)\ar[ur]^{(\FF_h)_*}\ar[rr]_{((\GG\circ\FF)_h)_*} & &
   \Coh_{(\GG\circ\FF)_h(h)}(C).}                                  \]
\end{cor}

\subsection{Homotopy equivalence of twisting elements}\label{hh'}

\begin{df}\label{df-hmtp}
Two twisting elements $h,h'$ in an $A_\infty$-algebra $(A,\{\bb_n\})$ are
{\em homotopic} (through an element $c$) if there exists $c\in A$ (of odd
degree) such that
\[  h'-h=\sum_{r,t\ge0}\bb_{r+1+t}(\ubrace{h'}{r},c,\ubrace{h}{t}). \]
\end{df}

Under this situation, define a map $\uppsi_c\colon A\to A$ by
\[  a\in A\mapsto\uppsi_c(a)=a+
     \sum_{r,t\ge0}\bb_{r+t+2}(\ubrace{h'}{r},c,\ubrace{h}{t},a). \]
We note that $\uppsi_c(h)=h'-\bd_{h'}c$.

\begin{lm}\label{tw-hmtp}
Suppose $A$ is an $A_\infty$-algebra and $h,h'\in A$ are two twisting elements.
Then 

(i) $h$ and $h'$ are homotopic through $c\in A$ if and only if
\[ \TT_{h'}(h')-\TT_h(h)=\bb(\TT_{h'}(1)\ox c\ox\TT_h(1)); \]

(ii) in this case, we have, for any $a\in A$,
\[ \TT_{h'}(\uppsi_c(a))-\TT_h(a)=
\bb(\TT_{h'}(c)\ox\TT_h(a))+\TT_{h'}(c)\ox\TT_h(\bd_ha).   \]
Moreover, $\uppsi_c(a)$ is the unique element satisfying this equality.
\end{lm}

\noindent{\em Proof:}
(i) Following the proof of Lemma~\ref{bdh} and using Definition~\ref{df-hmtp},
we get 
\begin{eqnarray*}
\vc\bb(\TT_{h'}(1)\ox c\ox\TT_h(1))                                   \\
 \eq\TT_{h'}(1)\ox\bd_{h'}h'\ox\TT_{h'}(1)\ox c\ox\TT_h(1)
 +\TT_{h'}(1)\ox(h'-h)\ox\TT_h(1)                                     \\
 \vc-\TT_{h'}(1)\ox c\ox\TT_{h'}(1)\ox\bd_hh\ox\TT_h(1)               \\
\eq\TT_{h'}(h')-\TT_h(h)
\end{eqnarray*}
since $\bd_{h'}h'=0=\bd_hh$.\\
(ii) follows from a similar calculation
\begin{eqnarray*}
\vc\bb(\TT_{h'}(c)\ox\TT_h(a))                                        \\
\eq\TT_{h'}(1)\ox\bd_{h'}h'\ox\TT_{h'}(1)\ox c\ox\TT_h(a)
 +\TT_{h'}(1)\ox(\uppsi_c(a)-a)+\TT_{h'}(1)\ox(h'-h)\ox\TT_h(1)\ox a  \\
\vc-\TT_{h'}(1)\ox c\ox\TT_h(1)\ox\bd_hh\ox\TT_h(a)
 -\TT_{h'}(c)\ox\TT_h(1)\ox\bd_ha                                     \\
\eq\TT_{h'}(\uppsi_c(a))-\TT_h(a)-\TT_{h'}(c)\ox\TT_h(\bd_ha). 
\end{eqnarray*}
The uniqueness is clear.
\qed

It can be shown that this homotopy is an equivalence relation \cite{Pr,K86}.
In fact, the equality in Lemma~\ref{tw-hmtp}(i) means that $\TT_h(h)$ and
$\TT_{h'}(h')$ are homotopic as coalgebra morphisms from $\kk$ to $\bart A$.
More precisely, we have

\begin{pro}\label{trans}
Let $A$ be an $A_\infty$-algebra.

(i) If $h,h'\in A$ are twisting elements that are homotopic through $c$,
then $\uppsi_c$ is a cochain map from $(A,\bd_h)$ to $(A,\bd_{h'})$, i.e.,
$\uppsi_c\circ\bd_h=\bd_{h'}\circ\uppsi_c$.

(ii) If $h''\in A$ is another twsiting element and $h',h''$ are homotopic 
through $c'$, then $h,h''$ are homotopic through
\[ c''=c+c'+\sum_{r,s,t\ge0}
       \bb_{r+s+t+2}(\ubrace{h''}{r},c',\ubrace{h'}{s},c,\ubrace{h}{t}).  \]

(iii) Under the above condition, $\uppsi_{c'}\circ\uppsi_c$ and 
$\uppsi_{c''}\colon(A,\bd_h)\to(A,\bd_{h''})$ are homotopic cochain maps.
\end{pro}

\noindent{\em Proof:}
(i) We use Lemma~\ref{tw-hmtp}(ii) in two ways.
Applying $\bb$ on the formula, we get
\[ \bb(\TT_{h'}(c)\ox\TT_h(\bd_ha))=
   \TT_{h'}(\bd_{h'}\uppsi_c(a))-\TT_h(\bd_ha). \]
On the other hand, replacing $a$ by $\bd_ha$ in the same formula, we get
\[\bb(\TT_{h'}(c)\ox\TT_h(\bd_ha))=\TT_{h'}(\uppsi_c(\bd_ha))-\TT_h(\bd_ha).\]
Therefore $\uppsi_c\circ\bd_h=\bd_{h'}\circ\uppsi_c$ on $A$.\\
(ii) First, we calculate
\begin{eqnarray*}
\vc\bb(\TT_{h''}(1)\ox c'\ox\TT_{h'}(1)\ox c\ox\TT_h(1))    \\
\eq\TT_{h''}(1)\ox(h''-h')\ox\TT_{h'}(1)\ox c\ox\TT_h(1)
   -\TT_{h''}(1)\ox c'\ox\TT_{h'}(1)\ox(h'-h)\ox\TT_h(1)    \\
\vc\quad+\TT_{h''}(1)\ox(c''-c-c')\ox\TT_h(1)               \\
\eq\TT_{h''}(1)\ox c''\ox\TT_h(1)-\TT_{h'}(1)\ox c\ox\TT_h(1)
   -\TT_{h''}(1)\ox c'\ox\TT_{h'}(1).
\end{eqnarray*}
Applying $\bb$ to both sides and using Lemma~\ref{tw-hmtp}(i), we get
\[ \bb(\TT_{h''}(1)\ox c''\ox\TT_h(1))=(\TT_{h'}(h')-\TT_h(h))
   +(\TT_{h''}(h'')-\TT_{h'}(h'))=\TT_{h''}(h'')-\TT_h(h).   \]
The result follows from Lemma~\ref{tw-hmtp}(i).\\
(iii) A similar calculation shows that, for any $a\in A$,
\begin{eqnarray*}
\vc\bb(\TT_{h''}(c')\ox\TT_{h'}(c)\ox\TT_h(a))
   +\TT_{h''}(c')\ox\TT_{h'}(c)\ox\TT_h(\bd_ha)            \\
\eq\TT_{h''}(c')\ox\TT_h(a)-\TT_{h'}(c)\ox\TT_h(a)
   -\TT_{h''}(c')\ox\TT_{h'}(\uppsi_c(a))+\TT_{h''}(\uppsi_{c',c}(a)),
\end{eqnarray*}
where the map $\uppsi_{c',c}\colon A\to A$ is defined by
\[ a\in A\mapsto\uppsi_{c',c}(a)=\sum_{r,s,t\ge0}
   \bb_{r+s+t+3}(\ubrace{h''}{r},c',\ubrace{h'}{s},c,\ubrace{h}{t},a).  \]
Applying $\bb$ to the above, we get, after simplifying and using (i), 
the desired relation
\[ \uppsi_{c'}(\uppsi_c(a))-\uppsi_{c''}(a)
=\bd_{h''}\uppsi_{c',c}(a)+\uppsi_{c',c}(\bd_ha).   \]
\qed

\begin{cor}\label{isom-c}
Let $A$ be an $A_\infty$-algebra and $h,h'\in A$, two twisting elements
that are homotopic through $c$.
Then

(i) $\uppsi_c\colon A\to A$ induces an isomorphism 
$(\uppsi_c)_*\colon \Coh_h(A)\to \Coh_{h'}(A)$;

(ii) if $h''$ is another twisting element and $h',h''$ are homotopic through
$c'$, then there is a commutative diagram
\[ \xymatrix@=1pc{& \Coh_{h'}(A)\ar[dr]^{(\uppsi_{c'})_*} &  \\
\Coh_h(A)\ar[ur]^{(\uppsi_c)_*}\ar[rr]_{(\uppsi_{c''})_*} & &
\Coh_{h''}(A),} \]
where $c''$ is given by Proposition~\ref{trans}(i).
\end{cor}

\noindent{\em Proof:}
As in \cite{K86}, $c'$ can be chosen such that $c''$ in
Proposition~\ref{trans}(ii) vanishes.
By Proposition~\ref{trans}(iii), $\uppsi_c$ induces an isomorphism 
$(\uppsi_c)_*\colon \Coh_h(A)\to \Coh_{h'}(A)$ with inverse $(\uppsi_{c'})_*$.
The rest follows.
\qed

\subsection{Morphisms on homotopic twisting elements}
Let $A$ be an $A_\infty$-algebra.
Suppose $h',h\in A$ satisfy the relation 
\[ h'-h=\sum_{r,t\ge0}\bb_{r+1+t}(\ubrace{\overline{h'}}{r},c,\ubrace{h}{t}) \]
for some $c\in A$ (cf.\ Definition~\ref{df-hmtp}, but $h,h'$ need not be
twisting elements).
Let $\ff\colon A\to B$ be a morphism of $A_\infty$-algebras.
We set
\[  \FF_{h',h}(c)=
    \sum_{r,t\ge0}\ff_{r+1+t}(\ubrace{\overline{h'}}{r},c,\ubrace{h}{t}) \]
and define a map $\FF_{h',h}(c,\cdot\,)\colon A\to A$ by
\[ a\in A\mapsto\FF_{h',h}(c,a)=\sum_{r,t\ge0}
                 \ff_{r+t+2}(\ubrace{\overline{h'}}{r},c,\ubrace{h}{t},a). \]

\begin{lm}\label{f-hmtp}
(i) In the above notations, we have the identities
\[ \FF(\TT^A_{\overline{h'}}(1)\ox c\ox\TT^A_h(1))
   =\TT^B_{\overline{\FF_{h'}(h')}}(1)\ox\FF_{h',h}(c)\ox\TT^B_{\FF_h(h)}(1).\]

(ii) For any $a\in A$, we have
\[ \FF(\TT^A_{\overline{h'}}(c)\ox\TT^A_h(a))
   =\TT^B_{\overline{\FF_{h'}(h')}}(\FF_{h',h}(c))\ox\TT^B_{\FF_h(h)}(\FF_h(a))
    +\TT^B_{\overline{\FF_{h'}(h')}}(\FF_{h',h}(c,a)).   \]
\end{lm}

\noindent{\em Proof:}
(i) Following the proof of Lemma~\ref{morh}, we get
\begin{eqnarray*}
\vc\FF(\TT^A_{\overline{h'}}(1)\ox c\ox\TT^A_h(1))                      \\
\eq\sum_{\substack{r\ge0;\,i_0,\dots,i_r>0\\ t\ge0;\,j_0\dots,j_t>0}}
  \ff_{i_r}(\ubrace{\overline{h'}}{i_r})\ox\cdots\ox\ff_{i_1}(\ubrace{h'}{i_1})
  \ox\ff_{i_0+j_0+1}(\ubrace{\overline{h'}}{i_0},c,\ubrace{h}{j_0})\ox  \\
\vc\quad\quad\quad\quad\quad\quad\quad\quad
   \ox\ff_{j_1}(\ubrace{h}{j_1})\ox\cdots\ox\ff_{j_t}(\ubrace{h}{j_t}) \\
\eq\sum_{r,t\ge0}\xbrace{\overline{\FF_{h'}(h')}}{r}\ox\,\FF_{h',h}(c)\ox
   \xbrace{\FF_h(h)}{t}                                                 \\   
\eq\TT^B_{\overline{\FF_{h'}(h')}}(1)\ox\FF_{h',h}(c)\ox\TT^B_{\FF_h(h)}(1).
\end{eqnarray*}
(ii) can be proved by a similar calculation.                \qed

\begin{thm}\label{fhh'}
Let $\ff\colon A\to B$ be a morphism of $A_\infty$-algebras.
Suppose $h,h'\in A$ are twisting elements that are homotopic through $c$.
Then

(i) $\FF_h(h),\FF_{h'}(h')\in B$ are twisting elements that are homotopic
through $\FF_{h',h}(c)$;

(ii) the diagram
\[ \xymatrix@=1pc{(A,\bd^A_h)\ar[dd]_{F_h}\ar[rrr]^{\uppsi^A_c}
                  & & & (A,\bd^A_{h'}) \ar[dd]^{F_{h'}} \\
   & & & \\
   (B,\bd^B_{\FF_h(h)})\ar[rrr]_{\uppsi^B_{FF_{h',h}(c)}} & & &
   (B,\bd^B_{\FF_{h'}(h')}).}     \]
commutes up to a cochain homotopy.
\end{thm}

\noindent{\em Proof:}
(i) Applying $\FF$ to the formula in Lemma~\ref{tw-hmtp}(i) and using
Lemma~\ref{f-hmtp}(i), we get
\begin{eqnarray*}
\vc\TT^B_{\FF_{h'}(h')}(\FF_{h'}(h'))-\TT^B_{\FF_h(h)}(\FF_h(h))
  =\FF(\bb^A(\TT^A_{h'}(1)\ox c\ox\TT^A_h(1)))                     \\
\eq\bb^B(\FF(\TT^A_{h'}(1)\ox c\ox\TT^A_h(1)))
=\bb^B(\TT^B_{\FF_{h'}(h')}(1)\ox\FF_{h',h}(c)\ox\TT^B_{\FF_h(h)}(1)).
\end{eqnarray*}
The result is then proved by using Lemma~\ref{tw-hmtp}(i) again.\\
(ii) Applying $\FF$ to the formula in Lemma~\ref{tw-hmtp}(ii) and using
Lemma~\ref{f-hmtp}(ii), we get
\begin{eqnarray*}
\vc\TT^B_{\FF_{h'}(h')}(\FF_{h'}(\uppsi^A_c(a)))-\TT^B_{\FF_h(h)}(\FF_h(a)) \\
\eq\bb^B(\TT^B_{\FF_{h'}(h')}(\FF_{h',h}(c))\ox\TT^B_{\FF_h(h)}(\FF_h(a))
    +\TT^B_{\FF_{h'}(h')}(\FF_{h',h}(c,a)))       \\
\vc\quad+\TT^B_{\FF_{h'}(h')}(\FF_{h',h}(c))\ox\TT^B_{\FF_h(h)}(\FF_h(\bd_ha))
   +\TT^B_{\FF_{h'}(h')}(\FF_{h',h}(c,\bd_ha))     \\
\eq\bb^B(\TT^B_{\FF_{h'}(h')}(\FF_{h',h}(c))\ox\TT^B_{\FF_h(h)}(\FF_h(a)))
   +\TT^B_{\FF_h(h)}(\bd_{\FF_h(h)}\FF_h(a))       \\ 
\vc\quad+\TT^B_{\FF_{h'}(h')}(\bd_{\FF_{h'}(h')}\FF_{h',h}(c,a)
   +\FF_{h',h}(c,\bd_ha)).     
\end{eqnarray*}
On the other hand, applying Lemma~\ref{tw-hmtp}(ii) to $B$, we get
\begin{eqnarray*}
\vc\TT^B_{\FF_{h'}(h')}(\uppsi^B_{\FF_{h',h}(c)}(\FF_h(a)))
 -\TT^B_{\FF_h(h)}(\FF_h(a))                                     \\
\eq\bb^B(\TT^B_{\FF_{h'}(h')}(\FF_{h',h}(c))\ox\TT^B_{\FF_h(h)}(\FF_h(a)))
 +\TT^B_{\FF_{h'}(h')}(\FF_{h',h}(c))\ox
   \TT^B_{\FF_h(h)}(\bd_{\FF_h(h)}\FF_h(a)).   
\end{eqnarray*}
Comparing the two calculations and by the injectivity of
$\TT^B_{\FF_{h'}(h')}$, we have
\[ \FF_{h'}(\uppsi^A_c(a))-\uppsi^B_{\FF_{h',h}(c)}(\FF_h(a))
   =\bd_{\FF_{h'}(h')}\FF_{h',h}(c,a)+\FF_{h',h}(c,\bd_ha),     \]
which establishes the desired homotopy through $\FF_{h',h}(c,\cdot\,)$.
\qed

\subsection{Homomorphisms on cohomology induced by homotopic morphisms}
Suppose $\{\ff_n\}$ and $\{\gg_n\}$ are two $A_\infty$-morphisms from 
$A$ to $B$ that are homotopic through $\{\hh_n\}$.
In terms of graded differential coalgebras, 
$\FF,\GG\colon\bart A\to\bart B$ are homotopic through $\HH$.
For $h\in A$, let $\HH_h\colon A\to B$ be a map defined by
\[  a\in A\mapsto\HH_h(a)=\sum_{n=0}^\infty\hh_{n+1}(\ubrace{h}{n},a).   \]

\begin{lm}\label{h-t}
(i) If $h\in A$, then
\[ \HH(\TT^A_h(h))=\TT^B_{\overline{\GG_h(h)}}(1)
                              \ox\HH_h(h)\ox\TT^B_{\FF_h(h)}(1). \]

(ii) In addition, for any $a\in A$,
\[ \HH(\TT^A_h(a))=\TT^B_{\overline{\GG_h(h)}}(\HH_h(a))+
  \TT^B_{\overline{\GG_h(h)}}(\HH_h(h))\ox\TT^B_{\FF_h(h)}(\FF_h(a)).   \]
\end{lm}

\noindent{\em Proof:}
(i) follows from a direct calculation
\begin{eqnarray*}
\vc\HH(\TT^A_h(h))       \\
\eq\sum_{\substack{r,t\ge0,s>0\\ i_1,\dots,i_r>0\\ j_1,\dots,j_t>0}}
   \gg_{i_1}(\ubrace{\bar h}{i_1})\ox\cdots\ox\gg_{i_r}(\ubrace{\bar h}{i_r})
   \ox\hh_s(\ubrace{h}{s})\ox
   \ff_{j_1}(\ubrace{h}{j_1})\ox\cdots\ox\ff_{j_t}(\ubrace{h}{j_t})          \\
\eq\sum_{r,t\ge0}\xbrace{\overline{\GG_h(h)}}{r}
    \ox\,\HH_h(h)\ox\xbrace{\FF_h(h)}{t} \\
\eq\TT^B_{\overline{\GG_h(h)}}(1)\ox\HH_h(h)\ox\TT^B_{\FF_h(h)}(1).
\end{eqnarray*}
The proof of (ii) is similar.   \qed

\begin{thm}\label{fgh}
Suppose $\{\ff_n\}$, $\{\gg_n\}$ are two $A_\infty$-morphisms from $A$ to $B$
that are homotopic through $\{\hh_n\}$.
Let $h$ be a twisting element in $A$. Then

(i) $\FF_h(h)$, $\GG_h(h)$ are twisting elements in $B$ that are homotopic
through $\HH_h(h)\in B$;

(ii) $\uppsi_{\HH_h(h)}\circ\FF_h$ and 
$\GG_h\colon(A,\bd_h)\to(B,\bd_{\GG_h(h)})$ are two maps of cochain complexes
that are homotopic.
\end{thm}

\noindent{\em Proof.}
(i) Applying $\GG-\FF=\bb^B\circ\HH+\HH\circ\bb^A$ to $\TT^A_h(h)\in\bart A$
and using Lemma~\ref{h-t}(i), we get, 
\begin{eqnarray*}
\TT^B_{\GG_h(h)}(\GG_h(h))-\TT^B_{\FF_h(h)}(\FF_h(h))
\eq\bb^B(\HH(\TT^A_h(h)))+\HH(\bb^A(\TT^A_h(h)))                \\
\eq\bb^B(\TT^B_{\GG_h(h)}(1)\ox\HH_h(h)\ox\TT^B_{\FF_h(h)}(1)).
\end{eqnarray*}
The result then follows from Lemma~\ref{tw-hmtp}(i).\\
(ii) Applying the same formula to $\TT^A_h(a)\in\bart A$ ($a\in A$) and using 
Lemma~\ref{h-t}(ii), we get
\begin{eqnarray*}
\vc\TT^B_{\GG_h(h)}(\GG_h(a))-\TT^B_{\FF_h(h)}(\FF_h(a))                 \\
\eq\bb^B(\TT^B_{\GG_h(h)}(\HH_h(a))+\TT^B_{\GG_h(h)}(\HH_h(h))\ox
   \TT^B_{\FF_h(h)}(\FF_h(a)))+\HH(\TT^A_h(\bd_ha))+\HH(\bb(\TT^A_h(h))) \\
\eq\bb^B(\TT^B_{\GG_h(h)}(\HH_h(h))\ox\TT^B_{\FF_h(h)}(\FF_h(a)))
   +\TT^B_{\GG_h(h)}(\HH_h(h))\ox\bb(\TT^B_{\FF_h(h)}(\FF_h(a)))         \\
\vc+\TT^B_{\GG_h(h)}(\bd_{\GG_h(h)}\HH_h(a)+\HH_h(\bd_ha)).              
\end{eqnarray*}
Therefore
\[ \GG_h-\uppsi_{\HH_h(h)}\circ\FF_h=
   \bd_{\GG_h(h)}\circ\HH_h+\HH_h\circ\bd_h   \]
by the uniqueness in Lemma~\ref{tw-hmtp}(ii).
\qed

\begin{cor}\label{fgh-isom}
Under the above conditions, there is a commutative diagram
\[ \xymatrix@=1pc{& \Coh_h(A)\ar[dl]_{(\FF_h)_*}\ar[dr]^{(\GG_h)_*} &  \\
\Coh_{\FF_h(h)}(B)\ar[rr]_{(\uppsi_{\HH_h(h)})_*} & & 
\Coh_{\GG_h(h)}(B).} \]
\end{cor}

\subsection{Twisting elements and twisted cohomology of $C_\infty$-algebras}
\label{C-twist}

When the $A_\infty$-structure is $C_\infty$, there are a number of 
simplifications.
We summarize the results in the following

\begin{pro}
Suppose $\kk$ is a field of characteristic $0$.

(i) If $A$ is a $C_\infty$-algebra, then any closed element $h\in A$ of even
degree is a twisting element.

(ii) If $\ff\colon A\to B$ is a morphism of $C_\infty$-algebras, then 
$\FF_h(h)=\ff_1(h)$.

(iii) If two $C_\infty$-algebra morphisms $\ff,\gg\colon A\to B$ are
homotopic through $\{\hh_n\}$, then $\HH_h(h)=\hh_1(h)$.
\end{pro}
\noindent{\em Proof:}
This is because $\bb_1(h)=0$ and
\[  \bb_n(\ubrace{h}{n})=\frac{1}{n}\,\bb_n([\xbrace{h}{n-1}]\shuffle[h])=0 \]
for all $n\ge2$.
Similarly, for all $n\ge2$, $\ff_n(\ubrace{h}{n})=0$ and
$\hh_n(\ubrace{h}{n})=0$.
Therefore $\FF_h(h)=\ff_1(h)$ and $\HH_h(h)=\hh_1(h)$.
\qed

\section{Higher Massey products on the cohomology of $A_\infty$-algebras}
\label{massey}

The triple Massey product \cite{Mas}, was generalized to the context of
$A_\infty$-algebras \cite{St70}.
We give an explicit construction of the higher Massey products, usually defined
for differential graded algebras \cite{Kr,May}, for $A_\infty$-algebras.
Furthermore, we introduce an equivalence relation on the set of defining
systems under which the Massey product takes the same value in the cohomology; 
this clarifies the dependency of the Massey product on the defining systems.
We establish some properties of the Massey products and the naturality under
morphisms of $A_\infty$-algebras. 
Finally, we study the triple and higher Massey products of $C_\infty$-algebras.

\subsection{Triple Massey product}\label{mas3}

We first review the definition of the triple Massey product on the cohomology
of an $A_\infty$-algebra $A$.
Given classes $\al_1,\al_2,\al_3\in \Coh(A)$, let 
$a_{01},a_{12},a_{23}\in A$ be their representatives, respectively.
Suppose $\bar\bb_2(\al_1,\al_2)=0=\bar\bb_2(\al_2,\al_3)$.
Then there are $a_{02},a_{13}\in A$ such that 
$\bb_2(a_{01},a_{12})=-\bb_1(a_{02})$ and 
$\bb_2(a_{12},a_{23})=-\bb_1(a_{13})$.
If $A$ is a differential graded algebra, then
$\bb_2(a_{01},a_{13})+\bb_2(a_{02},a_{23})$ would be a cocycle representing
the usual triple Massey product.
However, when $A$ is a general $A_\infty$-algebra, this expression is no
longer closed.
Instead, with a correction term from $\bb_3$, we define \cite{St70}
\[ \mu(a_{01},a_{12},a_{23};a_{02},a_{13})=
\bb_2(a_{01},a_{13})+\bb_2(a_{02},a_{23})+\bb_3(a_{01},a_{12},a_{23}). \]
It is closed since
\begin{eqnarray*}
&& \bb_1(\mu(a_{01},a_{12},a_{23};a_{02},a_{13}))          \\
\eq \bb_2(\bb_2(a_{01},a_{12}),a_{23})
    +\bb_2(\overline{a_{01}},\bb_2(a_{12},a_{23}))
    +\bb_1\bb_3(a_{01},a_{12},a_{23})                      \\
\eq 0.
\end{eqnarray*}
Therefore the generalization of the triple Massey product for
$A_\infty$-algebras should be defined by the cocycle
$\mu(a_{01},a_{12},a_{23};a_{02},a_{13})$.

We now study how it depends on the various choices made.
First, $a_{01},a_{12},a_{23}$ can each differ by a coboundary.
Suppose, for example, $a'_{12}=a_{12}+\bb_1(c_{12})$ is used instead, 
then we can choose $a'_{02}=a_{02}+\bb_2(a_{01},c_{12})$,
$a'_{13}=a_{13}+\bb_2(c_{12},a_{23})$ and, accordingly, the difference
in the cocycles is
\[ \mu(a_{01},a'_{12},a_{23};a'_{02},a'_{13})
  -\mu(a_{01},a_{12},a_{23};a_{02},a_{13})
  =\bb_1\bb_3(a_{01},c_{12},a_{23}). \]
So the two cocycles represent the same class in $\Coh(A)$.
In addition, each of $a_{02},a_{13}$ can differ by a cocycle; this results 
a difference in $\bar\bb_2(\al_1,\Coh(A))+\bar\bb_2(\Coh(A),\al_3)$
of the class represented by $\mu(a_{01},a_{12},a_{23};a_{02},a_{13})$.
We define the triple Massey product $\la\al_1,\al_2,\al_3\ra$ as the set of
the cohomology classes which arises from all such choices of $a_{02},a_{13}$.
Thus $\la\al_1,\al_2,\al_3\ra$ is an element of 
$\Coh(A)/(\bar\bb_2(\al_1,\Coh(A))+\bar\bb_2(\Coh(A),\al_3))$.

Finally, we establish the naturality of the triple Massey product.
Let $\ff=\{\ff_n\}\colon A\to B$ be a morphism of $A_\infty$-algebras.
If $\al_1,\al_2,\al_3\in \Coh(A)$ satisfy
$\bar\bb^A_2(\al_1,\al_2)=0=\bar\bb^A_2(\al_2,\al_3)$, then $\beta_1=
(\ff_1)_*\al_1,\beta_2=(\ff_1)_*\al_2,\beta_3=(\ff_1)_*\al_3\in\Coh(B)$
satisfy the corresponding relations 
$\bar\bb^B_2(\beta_1,\beta_2)=0=\bar\bb^B_2(\beta_2,\beta_3)$.
The latters are represented by $b_{01}=\ff_1(a_{01}),
b_{12}=\ff_1(a_{12}),b_{23}=\ff_1(a_{23})\in B$, respectively.
We can choose $b_{02}=\ff_1(a_{02})+\ff_2(a_{01},a_{12})$ and
$b_{13}=\ff_1(a_{13})+\ff_2(a_{12},a_{23})$.
A straightforward calculation shows that
\[ \mu^B(b_{01},b_{12},b_{23};b_{02},b_{13})=
\ff_1(\mu^A(a_{01},a_{12},a_{23};a_{02},a_{13})). \]
Thus $\la\beta_1,\beta_2,\beta_3\ra\supset(\ff_1)_*\la\al_1,\al_2,\al_3\ra$.

\subsection{Matric $A_\infty$-algebras}

If $A$ is any vector space over $\kk$, we denote by $\mat{m}(A)$ the
set of $m\times m$ matrices with components in $A$.
We write $\aaa=(a_{ij})_{1\le i,j\le m}$, where $a_{ij}\in A$ is the
$(i,j)$-component of $\aaa\in\mat{m}(A)$.
Given two vector spaces $A$ and $B$, we define a product
$\odot\colon\mat{m}(A)\ox\mat{m}(B)\to\mat{m}(A\ox B)$ by
\[ (\aaa\odot\bbb)_{ij}=\sum_{k=1}^m a_{ik}\ox b_{kj},\quad1\le i,j\le m,  \]
where $\aaa\in\mat{m}(A)$ and $\bbb\in\mat{m}(B)$.
This product is associative under the natural identification of tensor products
of vector spaces.

If $A$ is an $A_\infty$-algebra, then $\mat{m}(A)$ has a grading given by
$\mat{m}(A)^p=\mat{m}(A^p)$ and there is a collection of $\kk$-multilinear
maps $\bb^{(m)}_n\colon\mat{m}(A)^{\ox n}\to\mat{m}(A)$ ($n\ge1$) defined by
\[ \bb^{(m)}_n(\aaa_1,\dots,\aaa_n)=\bb_n(\aaa_1\odot\cdots\odot\aaa_n),  \]
where $\aaa_1,\dots,\aaa_n\in\mat{m}(A)$ and $\bb_n$ acts on 
$\mat{m}(A^{\ox n})$ component-wise.

\begin{pro}\label{mat-b}
If $(A,\{\bb_n\})$ is an $A_\infty$-algebra, then so is
$(\mat{m}(A),\{\bb^{(m)}_n\})$.
\end{pro}

\noindent{\em Proof:}
For any $\aaa_1,\dots,\aaa_n\in\mat{m}(A)$, we have
\begin{eqnarray*}
\vc\sum_{\substack{r,t\ge0,\,s\ge1\\r+s+t=n}}
   \bb^{(m)}_{r+1+t}(\overline{\aaa_1},\dots,\overline{\aaa_r},
   \bb^{(m)}_s(\aaa_{r+1},\dots,\aaa_{r+s}),\aaa_{r+s+1},\dots,\aaa_n)  \\
\eq\sum_{\substack{r,t\ge0,\,s\ge1\\r+s+t=n}}
   \bb_{r+1+t}(\overline{\aaa_1}\odot\cdots\odot\overline{\aaa_r}\odot
   \bb_s(\aaa_{r+1}\odot\cdots\odot\aaa_{r+s})\odot\aaa_{r+s+1}
   \odot\cdots\odot\aaa_n),
\end{eqnarray*}
which is zero by Definition~\ref{df-b} since $\bb_n$ acts on $\mat{m}(A)$
component-wise.
\qed

\begin{df}
If $(A,\{\bb_n\})$ is an $A_\infty$-algebra, the $m\times m$
{\em matric $A_\infty$-algebra} on $A$ is $(\mat{m}(A),\{\bb^{(m)}_n\})$.
\end{df}

More generally, if $R$ is a ring that contains the field $\kk$, then we have
an $A_\infty$-algebra $R\ox_\kk A$ that also has an $R$-module structure.
The maps $\bb^R_n\colon(R\ox_\kk A)^{\ox n}\to R\ox_\kk A$ is given by
\[ \bb^R_n(r_1\ox a_1,\dots,r_n\ox a_n)
      =(r_1\cdots r_n)\ox\bb_n(a_1,\dots,a_n), \]
where $r_1,\dots,r_n\in R$ and $a_1,\dots,a_n\in A$.
The matric $A_\infty$-algebra is the special case when $R=\mat{m}(\kk)$.

If $\ff=\{\ff_n\}\colon A\to B$ is a morphism of $A_\infty$-algebras, then 
we define $\kk$-multilinear maps 
$\ff^{(m)}_n\colon(\mat{m}(A))^{\ox n}\to\mat{m}(B)$ for all $n\ge1$ by
\[ \ff^{(m)}_n(\aaa_1,\dots,\aaa_n)=\ff_n(\aaa_1\odot\cdots\odot\aaa_n),  \]
where $\aaa_1,\dots,\aaa_n\in\mat{m}(A)$ and $\ff_n$ acts on
$\mat{m}(A^{\ox n})$ component-wise.

\begin{pro}\label{mat-f}
If $\ff=\{\ff_n\}\colon A\to B$ is a morphism of $A_\infty$-algebras,
then so is $\ff^{(m)}=\{\ff^{(m)}_n\}\colon\mat{m}(A)\to\mat{m}(B)$.
\end{pro}

\noindent{\em Proof:}
The proof is similar to that of Proposition~\ref{mat-b}, using
Definition~\ref{df-f}.
\qed

If $\ff,\gg\colon A\to B$ are two morphisms of $A_\infty$-algebras that are
homotopic through $\hh=\{\hh_n\}$, then we can define $\kk$-multilinear
maps $\hh^{(m)}_n\colon(\mat{m}(A))^{\ox n}\to\mat{m}(B)$ for all $n\ge1$ by
\[ \hh^{(m)}_n(\aaa_1,\dots,\aaa_n)=\hh_n(\aaa_1\odot\cdots\odot\aaa_n),  \]
where $\aaa_1,\dots,\aaa_n\in\mat{m}(A)$ and $\hh_n$ acts on 
$\mat{m}(A^{\ox n})$ component-wise.

\begin{pro}
If $\ff,\gg\colon A\to B$ are two morphisms of $A_\infty$-algebras that are
homotopic through $\hh=\{\hh_n\}$, then 
$\ff^{(m)},\gg^{(m)}\colon\mat{m}(A)\to\mat{m}(B)$ are homotopic through
$\hh^{(m)}=\{\hh^{(m)}_n\}$.
\end{pro}

\noindent{\em Proof:}
The proof is similar to that of Proposition~\ref{mat-b}, using
Definition~\ref{df-h}.
\qed

Let $\mat{m}_+(A)\subset\mat{m}(A)$ be the subspace of strictly 
upper-triangular matrices in $A$, that is, $\aaa=(a_{ij})\in\mat{m}_+(A)$
if $a_{ij}=0$ for all $i\ge j$.
Since $\mat{m}_+(A)$ is preserved by the product $\odot$, we have
 
\begin{cor}
(i) If $(A,\{\bb_n\})$ is an $A_\infty$-algebra, then so is 
$(\mat{m}_+(A),\{\bb^{(m)}_n\})$.

(ii) If $\ff=\{\ff_n\}\colon A\to B$ is a morphism of $A_\infty$-algebras,
then so is
\[ \ff^{(m)}=\{\ff^{(m)}_n\}\colon\mat{m}_+(A)\to\mat{m}_+(B). \]

(iii) If $\ff,\gg\colon A\to B$ are two morphisms of $A_\infty$-algebras
that are homotopic through $\hh=\{\hh_n\}$, then 
$\ff^{(m)},\gg^{(m)}\colon\mat{m}_+(A)\to\mat{m}_+(B)$ are homotopic through
$\hh^{(m)}=\{\hh^{(m)}_n\}$.
\end{cor}

If $A$ is a vector space over $\kk$, then we regard elements of
$A^{\op m}=\pbrace{A}{m}$ as column vectors and we write
$\xxx=(x_i)_{1\le i\le m}$, where $x_i\in A$.
If $A$, $B$ are two vector spaces, there is a multiplication
$\odot\colon\mat{m}(A)\ox B^{\op m}\to(A\ox B)^{\op m}$
from the usual matrix multiplication on vectors.
If $\aaa=(a_{ij})\in\mat{m}(A)$ and $\xxx=(x_i)_{1\le i\le m}\in B^{\op m}$,
then $\aaa\odot\xxx=(\sum_{j=1}^ma_{ij}\ox x_j)_{1\le i\le m}$.
This multiplication satisfies the standard associativity upon natural
identification of tensor products of vector spaces.
If $A$ is an $A_\infty$-algebra, then $A^{\op m}$ is an $A_\infty$-module
over $\mat{m}(A)$ or $\mat{m}_+(A)$.
If $\aaa\in\mat{m}(A)$ and $\xxx\in A^{\op m}$, we write
\[ \bd_\aaa\xxx=\sum_{n=0}^\infty\bb_n(\dbrace{\aaa}{n}\odot\,\xxx), \]
where $\bb_n$ acts on the column vectors component-wise.

\begin{lm}\label{mat-vec}
If $\aaa\in\mat{m}(A)$ and $\xxx\in A^{\op m}$, let 
$\tilde\aaa={\aaa\quad\xxx\choose0\quad0}\in\mat{m+1}(A)$.
Then
$\bd_{\tilde\aaa}\tilde\aaa={\bd_\aaa\aaa\;\;\bd_\aaa\xxx\choose0\;\;\quad0}$.
\end{lm}

\noindent{\em Proof:}
This follows from
$\overbrace{\tilde\aaa\odot\cdots\odot\tilde\aaa}^{n\;\text{times}}
  ={\overbrace{\scriptstyle{\aaa\odot\cdots\odot\aaa}}^{n\;\text{times}}
  \quad\overbrace{\scriptstyle{\aaa\odot\cdots\odot\aaa}}^{n-1\;\text{times}}
  \!\odot\xxx\choose0\quad\quad\quad\quad0\;\;}$ for any $n\ge1$.
\qed

Finally, using the inclusions
\[ \cdots\subset\mat{m}(A)\subset\mat{m+1}(A)\subset\cdots\quad\mbox{and}\quad
   \cdots\subset\mat{m}_+(A)\subset\mat{m+1}_+(A)\subset\cdots,    \]
we get the $A_\infty$-algebras $(\mat{\infty}(A),\{\bb^{(\infty)}_n\})$ and
$(\mat{\infty}_+(A),\{\bb^{(\infty)}_n\})$ as direct limits.
The results in the section hold also for $m=\infty$.

\subsection{Higher Massey products}

We generalize the triple Massey product for $A_\infty$-algebras \cite{St70}
discussed in \S\ref{mas3} to higher Massey products of $m$ cohomology classes.
We now use the labeling $0\le i,j\le m$ for the components of
$\aaa=(a_{ij})\in\mat{m+1}(A)$.
If $\aaa,\aaa'\in\mat{m+1}_+(A)$, we write $\aaa\approx\aaa'$ if
$a_{ij}=a'_{ij}$ for all $i,j$ except $i=0,j=m$.
We have a simple

\begin{lm}\label{approx}
If $\aaa,\aaa'\in\mat{m+1}_+(A)$ and $\bbb,\bbb'\in\mat{m+1}_+(B)$ satisfy
$\aaa\approx\aaa'$ and $\bbb\approx\bbb'$, then
$\aaa\odot\bbb=\aaa'\odot\bbb'$.
\end{lm}

\begin{df}
Let $(A,\{\bb_n\})$ be an $A_\infty$-algebra and $m\ge 2$.
The matrix $\aaa=(a_{ij})_{0\le i,j\le m}\in\mat{m+1}_+(A)$ is a
{\em defining system} for $\al_1,\dots,\al_m\in \Coh(A)$ if

(i) $\bd_\aaa\aaa\approx0$; let $\mu(\aaa)=(\bd_\aaa\aaa)_{0m}\in A$;

(ii) each $\al_i$ ($i=1,\dots,m$) is represented by $a_{i-1,i}\in A$,
which is closed by (i).
\end{df}

When $m=2$, a defining system of $\al_1,\al_2\in\Coh(A)$ is simply a choice
of representatives $a_{01},a_{12}\in A$ of $\al_1,\al_2$, respectively, and
$\mu(a_{01},a_{12})=\bb_2(a_{01},a_{12})$ is a cocycle representing the class
$\bar\bb_2(a_{01},a_{12})\in\Coh(A)$.
The case $m=3$ is about the triple Massey product discussed in \S\ref{mas3}.
For a general $m$, the formula is
\[ \mu(\aaa)=\sum_{r=1}^m\;\sum_{0=i_0<i_1<\cdots<i_r=m}
             \bb_r(a_{i_0i_1},a_{i_1i_2},\cdots,a_{i_{r-1}i_r}).   \]

We remark that if $\aaa\in\mat{m+1}_+(A)$, then
\[  \bd_\aaa\aaa=\sum_{n=1}^\infty\bb^{(m+1)}_n(\ubrace{\aaa}{n})
   \in\mat{m+1}_+(A)  \]
can be regarded as the ``curvature'' of $\aaa$ in the context of
(non-associative) $A_\infty$-algebras. 
Furthermore, $\mu(\aaa)=0$ if and only if $\aaa$ is a twisting element of 
$\mat{m+1}_+(A)$.
In this case, the equations $\bd_\aaa\aaa=0$ generalizes the Maurer-Cartan
equation for flat connections.
In general, applying Corollary~\ref{twisting}(i) to $\mat{m+1}_+(A)$, we get
\[ \sum_{r,t\ge0}\bb_{r+1+t}(\dbrace{\bar\aaa}{r}\odot\,
        \bd_\aaa\aaa\odot\dbrace{\aaa}{t})=0.                     \]
This is the non-associative, $A_\infty$-algebraic version of the Bianchi
identity.
(See \cite{Kr,May2,BT} for the case of differential graded algebras).
Furthermore, by Corollary~\ref{twisting}(iii), we get
\[ \bd_{\bar\aaa}\bd_\aaa\ccc=-\sum_{r,t\ge0}\bb_{r+t+2}
  (\dbrace{\bar\aaa}{r}\odot\,\bd_\aaa\aaa\odot\dbrace{\aaa}{t}\odot\,\ccc), \]
where $\ccc$ can be either in $\mat{m+1}(A)_+$ or in $A^{\op(m+1)}$.
This reflects the familiar relation in geometry between the square of the
connection and the curvature.

\begin{pro}\label{mu0m}
(i) If $\aaa$ is a defining system of $\al_1,\dots,\al_m\in \Coh(A)$,
then $\bb_1(\mu(\aaa))=0$.

(ii) If $\aaa'$ is another defining system and $\aaa'\approx\aaa$,
then $\mu(\aaa')-\mu(\aaa)=\bb_1(a'_{0m}-a_{0m})$
and hence $[\mu(\aaa')]=[\mu(\aaa)]\in \Coh(A)$.
\end{pro}

\noindent{\em Proof:}
(i) Since $\bd_\aaa\aaa\approx0$, by Lemma~\ref{approx}, the Bianchi identity
reduces to $\bb_1(\bd_\aaa\aaa)=0$ and hence $\bb_1(\mu(\aaa))=0$.\\
(ii) is obvious.

\begin{df}
The (higher) Massey product of the cohomology classes
$\al_1,\dots,\al_m\in \Coh(A)$ is {\em defined} if there is a defining system
$\aaa$ of them.
The class $[\mu(\aaa)]\in H(A)$ is the {\em (higher) Massey product} of
$\al_1,\dots,\al_m$ {\em through the defining system $\aaa$}.
The {\em (higher) Massey product} of $\al_1,\dots,\al_m\in \Coh(A)$
is the set $\la\al_1,\dots,\al_m\ra\subset \Coh(A)$ of elements
$[\mu(\aaa)]$, where $\aaa$ runs over all the defining systems of
$\al_1,\dots,\al_m$.
\end{df}

We establish a property of the higher Massey product, generalizing
the case when $A$ is a differential graded algebra \cite{Kr,May2}.

\begin{thm}
Let $A$ be an $A_\infty$-algebra.
Suppose the Massey product of $\al_1,\dots,\al_m\in\Coh(A)$ is defined.
Then for any $\gam\in\Coh(A)$, the Massey product of
$\overline{\al_1},\cdots,\overline{\al_{m-1}},\bar\bb_2(\al_m,\gam)$ is also
defined, and $\la\overline{\al_1},\cdots,\overline{\al_{m-1}},
\bar\bb_2(\al_m,\gam)\ra\supset-\bb_2(\la\al_1,\cdots,\al_m\ra,\gam)$.
\end{thm}

\noindent{\em Proof:}
Let $\aaa\in\mat{m+1}_+(A)$ be a defining system of $\al_1,\dots,\al_m$.
We write $\aaa={\aaa_0\quad\!\aaa'\choose0\quad0}$, where
$\aaa_0\in\mat{m}_+(A)$ and $\aaa'\in A^{\op m}$.
By Lemma~\ref{mat-vec}, $\bd_{\aaa_0}\aaa_0=0$ in $\mat{m}_+(A)$
and $\bd_{\aaa_0}\aaa'={\mu(\aaa)\choose\ooo}$ in $A^{\op m}$.
Let $c\in\ker\bb_1\subset A$ be a representative of $\gam$ and set
$\ccc={\ooo\choose c}\in A^{\op(m+1)}$.
Then $\bd_\aaa\ccc={\bbb'\choose 0}$ for some $\bbb'\in A^{\op m}$.
Let $\bbb={\overline{\aaa_0}\quad\!\bbb'\choose0\quad0}\in\mat{m+1}_+(A)$.
Then $b_{i-1,i}=\overline{a_{i-1,i}}$ for $i=1,\dots,m-1$ and
$b_{m-1,m}=\bb_2(a_{m-1,m},c)$, which are representatives of
$\overline{\al_1},\dots,\overline{\al_{m-1}},\bar\bb_2(\al_m,\gam)\in\Coh(A)$,
respectively.
Using Lemma~\ref{mat-vec} again, we get
\[ \bd_\bbb\bbb={\bd_{\overline{\aaa_0}}\overline{\aaa_0}\quad\!
                 \bd_{\overline{\aaa_0}}\bbb'\choose0\;\quad\quad0}
               =\left(\ooo\quad\!\bd_{\bar\aaa}\bd_\aaa\ccc\right)   \]
since $0=\overline{\bd_{\aaa_0}\aaa_0}=
-\bd_{\overline{\aaa_0}}\overline{\aaa_0}$.
By Corollary~\ref{twisting}(iii) and Lemma~\ref{approx}, we have
\[  \bd_{\bar\aaa}\bd_\aaa\ccc=-\bar\bb_2(\bd_\aaa\aaa\odot\ccc)
                  =-{\bb_2(\mu(\aaa),c)\choose\ooo}\in A^{\op(m+1)}.   \]
So $\bbb$ is a defining system of
$\overline{\al_1},\dots,\overline{\al_{m-1}},\bar\bb_2(\al_m,\gam)$
and $\mu(\bbb)=-\bb_2(\mu(\aaa),c)$.
\qed

A similar statement can be made for the Massey product of
$\bar\bb_2(\bar\gam,\al_1),\al_2,\dots,\al_m$.

\subsection{Homotopy equivalence of defining systems}

In Proposition~\ref{mu0m}, we saw that the cohomology class of $\mu(\aaa)$
does not depend on the $(0,m)$-component of $\aaa$.
In fact, it is invariant under a wider equivalence of the defining systems.

\begin{df}\label{aaa'}
Let $A$ be an $A_\infty$-algebra.
Two defining systems $\aaa,\aaa'$ of $\al_1,\dots,\al_m\in \Coh(A)$ are
{\em homotopic} (through $\ccc$) if there is $\ccc\in\mat{m+1}_+(A)$ such that
\[ \aaa'-\aaa=\sum_{r,t\ge0}\bb_{r+1+t}(\dbrace{\overline{\aaa'}}{r}
                                \odot\,\ccc\odot\dbrace{\aaa}{t}).     \]
$\aaa$ and $\aaa'$ are {\em equivalent} if there is $\ccc\in\mat{m+1}_+(A)$ 
such that
\[ \aaa'-\aaa\approx\sum_{r,t\ge0}\bb_{r+1+t}(\dbrace{\overline{\aaa'}}{r}
                                \odot\,\ccc\odot\dbrace{\aaa}{t}).     \]
\end{df}

Clearly, if $\aaa\approx\aaa'$ or if $\aaa$ and $\aaa'$ are homotopic, then
they are equivalent.
In fact, this equivalence is the weakest relation with this property. 
When $m=3$, the equivalence of the defining systems reflects the ambiguity
in defining the triple Massey product in \S\ref{mas3}. 

If $\aaa\in\mat{m+1}_+(A)$, then the map
$\TT^{(m+1)}_\aaa\colon\mat{m+1}_+(A)\to\bart\mat{m+1}_+(A)$ was
given in \S\ref{sec-tw}.
We define another map $\TT_\aaa\colon\mat{m+1}_+(A)\to\mat{m+1}_+(\bart A)$
by the composition of $\TT^{(m+1)}_\aaa$ with the product $\odot$.
That is, we have
\[ \TT_\aaa\colon\bbb\in\mat{m+1}_+(A)\mapsto
   \sum_{n=0}^\infty\dbrace{\aaa}{n}\odot\,\bbb. \]
Let $\TT_\aaa(1)=1+\TT_\aaa(\aaa)\in\mat{m+1}_+(\kk\op\bart A)$.
We recall that $\bb$ acts on $\mat{m+1}_+(\bart A)$ component-wise.

\begin{lm}\label{lm-a}
Let $(A,\{\bb_n\})$ be an $A_\infty$-algebra.

(i) If $\aaa\in\mat{m+1}_+(A)$ is a defining system, then
$\bb\circ\TT_\aaa=\TT_{\bar\aaa}\circ\bd_\aaa$ on $\mat{m+1}_+(A)$.

(ii) Two defining systems $\aaa,\aaa'\in\mat{m+1}_+(A)$ are homotopic through
$\ccc$ if and only if
\[ \TT_{\aaa'}(\aaa')-\TT_\aaa(\aaa)=
   \bb(\TT_{\overline{\aaa'}}(1)\odot\ccc\odot\TT_\aaa(1)).   \]
\end{lm}

\noindent{\em Proof:}
(i) By Lemma~\ref{bdh}(i), we have
\[ \bb^{(m+1)}\circ\TT^{(m+1)}_\aaa=\TT^{(m+1)}_{\bar\aaa}\circ\bd_\aaa
   +\TT^{(m+1)}_{\bar\aaa}(\bd_\aaa\aaa)\circ\TT^{(m+1)}_\aaa.   \]
Taking the product $\odot$, the second term on the right hand side vanishes 
by Lemma~\ref{approx} and we get the result.\\
(ii) Following the proof of Lemma~\ref{tw-hmtp}(i), we get
\begin{eqnarray*}
\vc\bb^{(m+1)}(\overline{\TT^{(m+1)}_{\aaa'}(1)}\ox\ccc
  \ox\TT^{(m+1)}_\aaa(1))                                               \\
\eq\TT^{(m+1)}_{\aaa'}(1)\ox\overline{\bd_{\aaa'}\aaa'}\ox
  \overline{\TT^{(m+1)}_{\aaa'}(1)}\ox\ccc\ox\TT^{(m+1)}_\aaa(1)
  +\TT^{(m+1)}_{\aaa'}(1)\ox(\aaa'-\aaa)\ox\TT^{(m+1)}_\aaa(1)          \\
 \vc-\TT^{(m+1)}_{\aaa'}(1)\ox\bar\ccc\ox\overline{\TT^{(m+1)}_\aaa(1)}
  \ox\bd_\aaa\aaa\ox\TT^{(m+1)}_\aaa(1).
\end{eqnarray*}
Taking the product $\odot$, the first and third terms on the right hand side
vanish by Lemma~\ref{approx} since
$\bd_\aaa\aaa\approx\bd_{\aaa'}\aaa'\approx0$.
The result follows.
\qed

\begin{thm}\label{aaa-mu}
If $\aaa$ and $\aaa'\in\mat{m+1}_+(A)$ are two equivalent defining systems,
then $ [\mu(\aaa')]=[\mu(\aaa)]\in \Coh(A)$.
\end{thm}

\noindent{\em Proof:}
By Proposition~\ref{mu0m}(ii), it suffices to show the result when $\aaa$ and
$\aaa'$ are homotopic.
By Lemma~\ref{lm-a}(i) and Lemma~\ref{approx}, we have
\[  \bb(\TT_\aaa(\aaa))=\TT_\aaa(\bd_\aaa\aaa)=\bd_\aaa\aaa.    \]
So, applying $\bb$ to the formula in Lemma~\ref{lm-a}(ii), we get
$\bd_{\aaa'}\aaa'=\bd_\aaa\aaa$ and $\mu(\aaa')=\mu(\aaa)$.
\qed

\subsection{Naturality of the higher Massey products}\label{nat-m}

Given a morphism $\ff=\{\ff_n\}\colon A\to B$ of $A_\infty$-algebras,
the induced morphism
$\ff^{(m+1)}=\{\ff^{(m+1)}_n\}\colon\mat{m+1}(A)\to\mat{m+1}(B)$
of $A_\infty$-algebras determines to a morphism
$\FF^{(m+1)}\colon\bart\mat{m+1}(A)\to\bart\mat{m+1}(B)$ of coalgebras.
For any $\aaa,\bbb\in\mat{m+1}(A)$, we have
\[ \FF^{(m+1)}_{\,\aaa}(\bbb)=
   \sum_{n=0}^\infty\ff_n(\dbrace{\aaa}{n}\odot\,\bbb),  \]
which we denote by $\FF_\aaa(\bbb)$ for short.

\begin{thm}\label{mas-f}
Suppose $\ff\colon A\to B$ of $A_\infty$-algebras and $m\ge2$.

(i) If $\aaa=(a_{ij})_{0\le i,j\le m}\in\mat{m+1}_+(A)$ is a defining system
of $\al_1,\dots,\al_m\in \Coh(A)$, then $\FF_\aaa(\aaa)\in\mat{m+1}_+(B)$
is a defining system of $(\ff_1)_*\al_1,\dots,(\ff_1)_*\al_m\in \Coh(B)$,
and $\mu^B(\FF_\aaa(\aaa))=\ff_1(\mu^A(\aaa))$.

(ii) If $\aaa$ and $\aaa'$ are two equivalent defining systems, then so
are $\FF_\aaa(\aaa)$ and $\FF_{\aaa'}(\aaa')$.
\end{thm}

\noindent{\em Proof:}
(i) Using Lemma~\ref{lm-a}(i) and following the proof of
Theorem~\ref{MORh}(ii), we have
\[ \bd_{\FF_\aaa(\aaa)}(\FF_\aaa(\aaa))=\FF_\aaa(\bd_\aaa\aaa).   \]
Since $\bd_\aaa\aaa\approx0$, the right hand side is equal to
$\ff_1(\bd_\aaa\aaa)\approx0$.
Therefore $\FF_\aaa(\aaa)$ is a defining system and
$\mu^B(\FF_\aaa(\aaa))=\ff_1(\mu^A(\aaa))$.
Next, we note that $(\FF_\aaa(\aaa))_{i-1,i}=\ff_1(a_{i-1,i})$.
Since $a_{i-1,i}$ is closed and represents $\al_i$ (for each $1\le i\le m$),
we get $\bb^B_1(\ff_1(a_{i-1,i}))=\ff_1(\bb^A_1(a_{i-1,i}))=0$ and that
$\ff_1(a_{i-1,i})$ represents $(\ff_1)_*\al_i$.\\
(ii) First, if $\aaa\approx\aaa'$, then 
$\FF_\aaa(\aaa)\approx\FF_{\aaa'}(\aaa')$.
If $\aaa$ and $\aaa'$ are homotopic through $\ccc$, then by the proofs of
Lemma~\ref{f-hmtp}(i) and Theorem~\ref{fhh'}(i)
(which can be adjusted without assuming twisting elements), we get
\begin{eqnarray*}
\vc\FF_{\aaa'}(\aaa')-\FF_\aaa(\aaa)        \\
\eq\sum_{r,t\ge0}\bb^B_{r+1+t}(
    \dbrace{\overline{\FF_{\aaa'}(\aaa')}}{r}\odot\,
    \FF_{\aaa',\aaa}(\ccc)\odot\dbrace{\FF_\aaa(\aaa)}{t}),
\end{eqnarray*}
where
\[ \FF_{\aaa',\aaa}(\ccc)=\sum_{r,t\ge0}\ff_{r+1+t}
       (\dbrace{\overline{\aaa'}}{r}\odot\,\ccc\odot\dbrace{\aaa}{t}).    \]
By Lemma~\ref{lm-a}(ii), $\FF_\aaa(\aaa)$ and $\FF_{\aaa'}(\aaa')$ are
homotopic through $\FF_{\aaa',\aaa}(\ccc)$.
\qed

\begin{cor}
Let $\ff\colon A\to B$ be a morphism of $A_\infty$-algebras.
If the Massey product of $\al_1\dots,\al_m\in\Coh(A)$ is defined, then
so is that of $(\ff_1)_*(\al_1)\dots,(\ff_1)_*(\al_m)\in\Coh(B)$, and
$\la(\ff_1)_*(\al_1)\dots,(\ff_1)_*(\al_m)\ra
\supset(\ff_1)_*\la\al_1\dots,\al_m\ra$.
\end{cor}

This generalizes the naturality of the triple Massey product in \S\ref{mas3}.

Recall that there is an $A_\infty$-structure $\{\bar\bb_n\}$ on $\Coh(A)$
after choosing the maps $\pp_1$ and $\qq_1$ (\S\ref{sec-ha}).
When $A$ is a differential graded algebra, it was a folklore that
the $A_\infty$-structure on $\Coh(A)$ gives the Massey products
\cite{St2,K80,Pr2,LPWZ}.
The precise statement seems to be a long standing puzzle.\footnote{It was
claimed in Theorem~3.1 of \cite{LPWZ} that when $A$ is a differential graded
algebra, the conclusion of Proposition~\ref{ha-mas} holds under a weaker
assumption that for any $j-i<m$, the Massey product of $\al_{i+1},\dots,\al_j$
is defined and contains $0$. We think that the stronger condition 
$\;\bar\bb_{j-i}(\al_{i+1},\dots,\al_j)=0$ as in Propistion~\ref{ha-mas}
is necessary even when $A$ is a differential graded algebra.}
We now establish the exact relationship in the more general context of 
$A_\infty$-algebras.

\begin{pro}\label{ha-mas}
Let $A$ be an $A_\infty$-algebra and $\al_1,\dots,\al_m\in\Coh(A)$.
If $\;\bar\bb_{j-i}(\al_{i+1},\dots,\al_j)=0$ for any $i,j$ satisfying
$0\le i<j\le m$, $j-i<m$, then the Massey product of $\al_1,\dots,\al_m$
is defined and $\bar\bb_m(\al_1,\dots,\al_m)\in\la\al_1,\dots,\al_m\ra$.
\end{pro}

\noindent{\em Proof:}
Notice that $\Coh(\Coh(A))=\Coh(A)$ as $\bar\bb_1=0$.
Let $\Al\in\mat{m+1}_+(\Coh(A))$ be given by $\al_{i-1,i}=\al_i$
($1\le i\le m$) and $\al_{ij}=0$ if $j\ne i+1$.
By the assumption, $\Al$ is a defining system of
$\al_1,\dots,\al_m\in\Coh(\Coh(A))=\Coh(A)$ and
$\mu^{\Coh(A)}(\Al)=\bar\bb_m(\al_1,\dots,\al_m)$.
It suffices to show that the Massey product of $\al_1,\dots,\al_m$ is defined
and contains $\mu^{\Coh(A)}(\Al)$.
Applying Theorem~\ref{mas-f}(i) to the quasi-isomorphism
$\qq\colon\Coh(A)\to A$ in \S\ref{sec-ha}, we get a defining system
$\aaa=\QQ_\Al(\Al)\in\mat{m+1}_+(A)$ of $[\qq_1(\al_i)]=\al_i$ ($1\le i\le m$).
More explicitly, $a_{ij}=\qq_{j-i}(\al_{i+1},\dots,\al_j)$ for $0\le i<j\le m$.
Therefore the Massey product of $\al_1,\dots,\al_m$ is defined.
Furthermore, $\qq_1(\mu^{\Coh(A)}(\Al))=\mu^A(\aaa)$ and hence
$\mu^{\Coh(A)}(\Al)=[\mu^A(\aaa)]\in\la\al_1,\dots,\al_m\ra$.
\qed

It is clear from the proof that the Massey product of $\al_1,\dots,\al_m$
as elements of $\Coh(\Coh(A))$ through the defining system $\Al$ is
identical to that as elements of $\Coh(A)$ through the defining system
$\aaa=\QQ_\Al(\Al)$.

Finally, we consider two morphisms $\ff,\gg\colon A\to B$ of
$A_\infty$-algebras that are homotopic through $\hh$.
Recall that $\ff^{(m+1)},\gg^{(m+1)}\colon\mat{m+1}_+A\to\mat{m+1}_+(B)$ are
morphisms of $A_\infty$-algebras that are homotopic through $\hh^{(m+1)}$.
As before, we use $\HH_\aaa(\aaa)$ to denote
\[ \HH^{(m+1)}_\aaa(\aaa)=\sum_{n=1}^\infty\hh_n(\dbrace{\aaa}{n}).  \]

\begin{thm}
Suppose $\ff,\gg\colon A\to B$ are two morphisms of $A_\infty$-algebras
that are homotopic through $\hh$.
Let $\aaa\in\mat{m+1}_+(A)$ be a defining system of
$\al_1\dots,\al_m\in \Coh(A)$, where $m\ge2$.
Then $\FF_\aaa(\aaa)$ and $\GG_\aaa(\aaa)$ are two equivalent defining
systems of $(\ff_1)_*\al_i=(\gg_1)_*\al_i\in\Coh(B)$ ($1\le i\le m$) and
hence define the same Massey product.
\end{thm}

\noindent{\em Proof:}
Using Lemma~\ref{h-t}(i) and following the proof of Theorem~\ref{fgh}(i),
we get
\begin{eqnarray*}
\vc\TT_{\GG_\aaa(\aaa)}(\GG_\aaa(\aaa))-\TT_{\FF_\aaa(\aaa)}(\FF_\aaa(\aaa))\\
\eq\bb^B(\TT^B_{\overline{\GG_\aaa(\aaa)}}\odot
    \HH_\aaa(\aaa)\odot\TT^B_{\FF_\aaa(\aaa)})+\HH(\TT^A_\aaa(\bd_\aaa\aaa)).
\end{eqnarray*}
Since $\HH(\TT^A_\aaa(\bd_\aaa\aaa))=\HH(\bd_\aaa\aaa)\approx0$,
the result follows from Lemma~\ref{lm-a}(ii) and Theorem~\ref{aaa-mu}.
\qed

\subsection{Massey products on the cohomology of $C_\infty$-algebras}

Let $(A,\{\bb_n\})$ be a $C_\infty$-algebra.
We consider first the triple Massey product (\S\ref{mas3}).
Let $\al_1,\al_2,\al_3\in\Coh(A)$ and let $a_{01},a_{12},a_{23}\in A$ be
their representatives, respectively.
Suppose
$\bar\bb_2(\al_1,\al_2)=\bar\bb_2(\al_2,\al_3)=\bar\bb_2(\al_3,\al_1)=0$.
Then there exist $a_{02},a_{13},a'_{21}\in A$ such that
$\bb_2(a_{01},a_{12})=-\bb_1(a_{02})$, $\bb_2(a_{12},a_{23})=-\bb_1(a_{13})$
and $\bb_2(a_{23},a_{01})=-\bb_1(a'_{21})$.
Recall that 
\[  [\al_1\ox\al_2]\shuffle\al_3=\al_1\ox\al_2\ox\al_3+(-1)^{|\al_2||\al_3|}
  \al_1\ox\al_3\ox\al_2+(-1)^{(|\al_1|+|\al_2|)|\al_3|}\al_3\ox\al_1\ox\al_2; \]
the same formula can be written for $[a_{01}\ox a_{12}]\shuffle a_{23}$.
We can check, using $\bb_2(A\shuffle A)=0$ and $\bb_3([A\ox A]\shuffle A)=0$,
that
\begin{eqnarray*}
\mu(a_{01},a_{12},a_{23};a_{02},a_{13})\!\!\!\!\!\!\!\!\vc
 +(-1)^{|\al_2||\al_3|}\mu(a_{01},a_{23},a_{12};
      (-1)^{|\al_1||\al_3|}a'_{21},(-1)^{|\al_2||\al_3|}a_{13})        \\
\vc+(-1)^{(|\al_1|+|\al_2|)|\al_3|}\mu(a_{23},a_{01},a_{12};a'_{21},a_{02})=0.
\end{eqnarray*}
Consequently, we have
\[ 0\in\la\al_1,\al_2,\al_3\ra+(-1)^{|\al_2||\al_3|}
\la\al_1,\al_3,\al_2\ra+(-1)^{(|\al_1|+|\al_2|)|\al_3|}\la\al_3,\al_1,\al_2\ra. \]
Similarly, corresponding to $\al_1\shuffle[\al_2\ox\al_3]$, we have
\[ 0\in\la\al_1,\al_2,\al_3\ra+(-1)^{|\al_1||\al_2|}
\la\al_2,\al_1,\al_3\ra+(-1)^{(|\al_2|+|\al_3|)|\al_1|}\la\al_2,\al_3,\al_1\ra. \]

We now consider higher Massey products.
Recall the notations $S_{r,n}$ and $\veps(\sig)$ in \S\ref{C-alg-mor}.

\begin{pro}
Let $A$ be a $C_\infty$-algebra and let $\al_1,\dots,\al_m\in\Coh(A)$.
Fix $r<m$.
Suppose for any $\sig\in S_{r,m}$ and any $i,j=0,\dots,m$ with $0<j-i<m$,
we have $\bar\bb_{j-i}(\al_{\sig(i+1)},\dots,\al_{\sig(j)})=0$.
Then the Massey product of
$\al_{\sig(1)},\dots,\al_{\sig(m)}$ is defined for any $\sig\in S_{r,m}$
and we have
\[ 0\in\sum_{\sig\in S_{r,m}}(-1)^{\veps(\sig)}\,
   \la\al_{\sig(1)},\dots,\al_{\sig(m)}\ra.  \]
\end{pro}

\noindent{\em Proof:}
By Proposition~\ref{ha-mas}, the assumption implies that for any 
$\sig\in S_{r,m}$, the Massey product of $\al_{\sig(1)},\dots,\al_{\sig(m)}$
is defined and $\bar\bb_m(\al_{\sig(1)},\dots,\al_{\sig(m)})\in
\la\al_{\sig(1)},\dots,\al_{\sig(m)}\ra$.
Since $(\Coh(A),\{\bar\bb_n\})$ is a $C_\infty$-algebra \cite{CG}, we have 
\[   0=\sum_{\sig\in S_{r,m}}(-1)^{\veps(\sig)}\,
       \bar\bb_m(\al_{\sig(1)},\dots,\al_{\sig(m)})  \]
and the result follows.
\qed

\section{Spectral sequence and higher Massey products}\label{spectral}

\subsection{A spectral sequence for the twisted cohomology}\label{spec3}

So far, the $A_\infty$-algebra $A$ is either $\Z$- or $\Z_2$-graded.
In both cases, we write $A=A^{\bar0}\op A^{\bar1}$, where
$A^{\bar0},A^{\bar1}$ are the even, odd parts of $A$, respectively.
Here $\bar k$ means the integer $k$ modulo $2$.
If $h\in A^{\bar0}$ is a twisting element, then the twisted differential
$\bd_h\colon A^{\bar k}\to A^{\overline{k+1}}$ defines a $\Z_2$-graded
cochain complex $(A^\bullet,\bd_h)$ and the twisted cohomology 
$\Coh_h(A)$ is always $\Z_2$-graded.

We now assume that $A$ is $\Z$-graded, i.e., $A=\bigoplus_{k\in\Z}A^k$.
If the degree of the twisting element $h$ is non-negative, then both the
component $h_0$ in $A^0$ and the positive-degree part $h-h_0$ are twisting
elements.
The twisted cohomology groups defined by $h_0$ and by $h-h_0$ behave very
differently (see \cite{MW} for the case of twisted de Rham
complex\footnote{In geometry, $\bd_h$ is a superconnection while $\bd_{h_0}$
is a usual connection, as $\sus^{-1}h_0$ is of degree $1$.}).
We now assume $h_0=0$, i.e., the twisting element $h$ has positive degree.
Denote the twisted differential $\bd_h$ on $A$ by $\dd$.
Then there is a natural filtration 
$\Fi$ of the $\Z_2$-graded cochain complex $(A^\bullet,\dd)$ given by
$\Fi^pA=\bigoplus_{n=p}^\infty A^n$, or
\[   \Fi^pA^{\bar k}=\bigoplus_{\substack{n\ge p\\ n=k\!\!\!\!\mod2}}A^n    \]
(see \cite{RW,AS,MW,LLW} when $A$ is the de Rham complex).
The graded components are $\Gr^pA=A^p$, or
\[  \Gr^pA^{\bar k}=\two{A^p}{p=k\!\!\!\!\mod 2}{0}{p\ne k\!\!\!\!\mod2.}   \]

\begin{lm}\label{E0}
There is a spectral sequence $\{E^{p\bar q}_r,\dd_r\}$ converging to
twisted cohomology $\Coh_h(A)$.
Moreover,

(i) $E^{p\bar 1}_r=0$ for any $p\in\Z$ and $r\ge0$;

(ii) $\dd_r=0$ if $r$ is even;

(iii) $E^{p\bar0}_{2m}=E^{p\bar0}_{2m+1}$ for any $p\in\Z$ and $m\ge1$.

(iv) $E^{p\bar0}_2=E^{p\bar0}_3=\Coh^p(A)$.
\end{lm}

\noindent{\em Proof:}
The filtration is clearly exhaustive and weakly convergent.
By Theorem~3.2 of \cite{Mc}, the corresponding spectral sequence converges
to the twisted cohomology group.\\
(i) $E^{p\bar 1}_r=0$ for all $r\ge0$ since 
$E^{p\bar 1}_0=\Gr^pA^{\overline{p+1}}=0$.\\
(ii) If $r$ is even, then either $q$ or $q-r+1$ is odd.
So $\dd_r\colon E^{p\bar q}_r\to E^{p+r,\overline{q-r+1}}_r$ is zero.\\
(iii) follows from (ii).\\
(iv) We have $E^{p\bar0}_1=E^{p\bar0}_0=A^p$ and
$\dd_1\colon A^p\to A^{p+1}$ is the (untwisted) differential $\bd=\bb_1$.
So $E^{p\bar0}_2=\Coh^p(A)$ and the rest follows from (iii).
\qed

We postpone the discussion on the general $\dd_r$ to \S\ref{d_r}.
Instead, generalizing the work on the twisted de Rham complex \cite{RW,AS},
we describe the differentials $\dd_3$ and $\dd_5$ when $A$ is an
$A_\infty$-algebra.
For simplicity, we assume that the twisting element $h$ is of homogeneous
degree $2$ (corresponding a closed $3$-form for the de Rham complex).
Then $\bb_n(\ubrace{h}{n})=0$ for any $n\ge1$.
Any element of $E^{p\bar0}_3$ is a class $[x]\in\Coh^p(A)$ represented by
a closed element $x\in A^p$, i.e., $\bb_1(x)=0$.
The map $\dd_3$ is given by $\dd_3[x]=[\bb_2(h,x)]$.
It is easy to check that $\bb_2(h,x)$ is closed (since $h$ and $x$ both are)
and the class $[\bb_2(h,x)]\in\Coh^{p+3}(A)$ does not depend on the choice
of the representative of $[x]$.

If $\dd_3[x]=[0]$, i.e., $\bb_2(h,x)=-\bb_1(x')$ for some $x'\in A^{p+2}$,
then $[x]$ represents a class $[x]_5\in E^{p\bar0}_4=E^{p\bar0}_5$.
When $A$ is a differential graded algebra, then $\dd_5[x]_5=[\bb_2(h,x')]_5$
(see \cite{RW,AS} for the case of de Rham complex).
But for a general $A_\infty$-algebra, we claim that
$\dd_5[x]_5=[\bb_2(h,x')+\bb_3(h,h,x)]_5$.
Note that $\bb_2(h,x')+\bb_3(h,h,x)=\mu(h,h,x;0,x')$ is closed and represents
a class in the triple Massey product $\la[h],[h],[x]\ra$.
Furthermore,
\[ \bb_2(h,\bb_2(h,x')+\bb_3(h,h,x))=-\bb_1(\bb_3(h,h,x')+\bb_4(h,h,h,x))  \]
is zero in $E^{p+7,\bar0}_3$ and hence
$[\bb_2(h,x')+\bb_3(h,h,x)]\in\Coh^{p+5}(A)=E^{p+5,\bar0}_3$ indeed
descends to a class in $E^{p+5,\bar0}_4=E^{p+5,\bar0}_5$.
We note that $x'$ is not unique; we can add an arbitrary closed element to
$x'$.
Yet the freedom in $x'$ does not exhaust all the defining systems of
$[h],[h],[x]$.
This results in a smaller ambiguity of the triple Massey product in
$\bar\bb_2([h],\Coh(A))$ (rather than
$\bar\bb_2([h],\Coh(A))+\bb_2(\Coh(A),[x])$ as in \S\ref{mas3}).
Happily, it is the subspace $\bar\bb_2([h],\Coh(A))$ that descends to zero
in $E_4=E_5$.

\subsection{Twisting elements and matric $A_\infty$-algebras}

We continue to assume that $A$ is a $\Z$-graded $A_\infty$-algebra and has
the filtration given in \S\ref{spec3}. 
If $a\in A$ is an element of positive, even degree, we write 
$a=a_2+a_4+a_6+\cdots$, where $a_{2i}\in A^{2i}$.
We define $\aaa^{(m)}=(a_{ij})_{0\le i,j\le m-1}\in\mat{m}_+(A)$ as the
strictly upper-triangular matrix given by $a_{ij}=a_{2(j-i)}$ if $i<j$.
Likewise, if $x\in\Fi^pA^{\bar p}$, we write $x=x_p+x_{p+2}+\cdots$, where
$x_{p+2i}\in A^{p+2i}$ ($i\ge0$).
We define a column vector
$\xxx_p^{(m)}=(x_{p-2(m-i-1)})_{0\le i\le m-1}\in A^{\op m}$.
In components, we have
\[ \aaa^{(m)}=\left(
\substack{0\quad a_2\quad a_4\quad\cdots\quad a_{2m-4}\;\;a_{2m-2} \\
\quad\;\;0\quad\;\;a_2\quad\cdots\quad a_{2m-6}\;\; a_{2m-4}\\
\quad\quad\;\;\ddots\;\;\;\ddots\quad\vdots\quad\quad\vdots\\
\quad\quad\quad\quad\quad\quad\;\;a_2\quad\;\;a_4\\
\quad\quad\quad\quad\quad\quad\;\;0\quad\;\;\;\,a_2\\
\quad\quad\quad\quad\quad\quad\quad\quad\;\;\;0\\}    \right),\quad
\xxx_p^{(m)}=\left(\substack{x_{p+2m-2} \\ x_{p+2m-4} \\ \vdots \\ x_{p+4} \\
   x_{p+2} \\ x_p}\right).   \]
We denote the zero matrix and the zero vector by $\ooo^{(m)}$ and
$\ooo_p^{(m)}$, respectively.

Given two graded $\Z$-graded vector spaces $A,B$ and two elements
$a\in A^{\bar0},\,b\in B^{\bar0}$ of positive degrees, then for each $m\ge2$,
the matrix that corresponds to $a\ox b\in(A\ox B)^{\bar0}$ is
$\aaa^{(m)}\odot\bbb^{(m)}$.
Similarly, if $x\in\Fi^pB^{\bar p}$, then for each $m\ge2$, $a\ox x$ determines
the column vectors $\aaa^{(m)}\odot\xxx_p^{(m)}\in(A\ox B)^{\op m}$.
Finally, by taking a direct limit, we can define $\aaa^{(\infty)}$ and
$\xxx_p^{(\infty)}$ with similar properties.

\begin{lm}\label{mat-bd}
Let $A$ be an $A_\infty$-algebra.
If $h\in A^{\bar0}$ is of positive degree and $x\in\Fi^pA^{\bar p}$,
$y\in\Fi^{p+1}A^{\overline{p+1}}$, then

(i) $\bd_hx-y\in\Fi^{p+2m+1}A$ if and only if 
$\bd_{\hhh^{(m)}}\xxx_p^{(m)}=\yyy_{p+1}^{(m)}$.
In particular, $\bd_hx=y$ if and only if 
$\bd_{\hhh^{(m)}}\xxx_p^{(m)}=\yyy_{p+1}^{(m)}$ for all $m\ge2$.

(ii) $\bd_hh\in\Fi^{p+2m+1}A$ if and only if
$\bd_{\hhh^{(m)}}\hhh^{(m)}=\ooo^{(m)}$.
In particular, $h$ is a twisting element in $A$ if and only if $\hhh^{(m)}$
is one in $\mat{m}_+(A)$ for any $m\ge2$.
\end{lm}

\noindent{\em Proof:}
(i) If $\bd_hx-y\in\Fi^{p+2m+1}A$, by comparing the homogeneous components
of $\bd_hx$ and $y$ in $A^{p+1}\op A^{p+3}\op\cdots\op A^{p+2m-1}$,
we get, for each $i=1,\dots,m$,
\[ \sum_{r=1}^i\sum_{\substack{i_1,\dots,i_r\ge1\\ i_1+\cdots+i_r=i}}
\bb_r(h_{2i_1},\dots,h_{2i_{r-1}},x_{p+2i_r})=0.  \]
These equalities are equivalent to
$\bd_{\hhh^{(m)}}\xxx_p^{(m)}=\yyy_{p+1}^{(m)}$. 
The rest is straightforward.\\
(ii) is proved similarly.
\qed

Combining Lemmas~\ref{mat-vec} and \ref{mat-bd}, we get

\begin{cor}\label{mat-df}
Suppose $h\in A^{\bar0}$ is of positive degree and $x\in\Fi^pA^{\bar p}$.
Then
${\hhh^{(m)}\;\;\xxx_p^{(m)}\choose\!\!0\quad\;\;\;0\;\;}\in\mat{m+1}_+(A)$
is a defining system if and only if $\bd_hh\in\Fi^{2m+1}A$ and
$\bd_hx\in\Fi^{p+2m-1}A$.
In this case, the component of $\bd_hx$ in $A^{p+2m-1}$ is
\[ (\bd_hx)_{p+2m-1}=
      \mu{\hhh^{(m)}\;\;\xxx_p^{(m)}\choose\!\!0\quad\;\;\;0\;\;}. \]
In particular, ${\hhh^{(m)}\;\;\xxx_p^{(m)}\choose\!\!0\quad\;\;\;0\;}$
is a twisting element if and only if $\bd_hh\in\Fi^{2m+1}A$ and
$\bd_hx\in\Fi^{p+2m+1}A$.
\end{cor}

\subsection{Higher differentials and higher Massey products}\label{d_r}

Let $A$ be an $A_\infty$-algebra and $h\in A$, a twisting element of positive
degree.
With the filtration $\Fi^pA$ defined in \S\ref{spec3}, we recall that
$\Gr^pA=A^p$ and the results in Lemma~\ref{E0}.
We apply the Appendix on spectral sequences to the $\Z_2$-graded setting of
twisted differential.
We want to describe the spaces 
$B^{p\bar0}_{2m+1}\subset Z^{p\bar0}_{2m+1}\subset A^p$,
$E^{p\bar0}_{2m+1}=Z^{p\bar0}_{2m+1}/B^{p\bar0}_{2m+1}$ and 
the maps $\dd_{2m+1}\colon E^{p\bar0}_{2m+1}\to E^{p+2m+1,\bar0}_{2m+1}$
for all $m\ge0$.

\begin{thm}\label{spec-thm}
Let $x\in A^p$.
Then

(i) $x_p\in Z^{p\bar0}_{2m+1}$ if and only if there exists 
$x=x_p+x_{p+2}+\cdots+x_{p+2m-2}$, where $x_{p+2i}\in A^{p+2i}$ ($0\le i<m$),
such that $\bd_{\hhh^{(m)}}\,\xxx_p^{(m)}=\ooo_{p+1}^{(m)}$, or equivalently, 
${\hhh^{(m)}\;\;\xxx_p^{(m)}\choose\!\!0\quad\;\;\;0\;\;}\in\mat{m+1}_+(A)$
is a twisting element.
In this case, $0\in\la\ubrace{[h_2]}{m-1},[x_p]\ra$.

(ii) $x_p\in B^{p\bar0}_{2m+1}$ if and only if there exists
$y=y_{p-2m+1}+y_{p-2m+3}+\cdots+y_{p-1}$, where $y_{p-2i+1}\in A^{p-2i+1}$
($0\le i<m$), such that $\bd_{\hhh^{(m)}}\,\yyy_{p-2m+1}^{(m)}=
\left(\substack{x_p\\ \ooo_{p-2m+2}^{(m-1)}}\right)$, or equivalently,
${\hhh^{(m)}\;\;\yyy_p^{(m)}\choose\!\!0\quad\;\;\;0\;\;}\in\mat{m+1}_+(A)$
is a defining system (of $\ubrace{[h_2]}{m-1},[y_{p-2m+1}]$) and 
$x_p=\mu{\hhh^{(m)}\;\;\yyy_p^{(m)}\choose\!\!0\quad\;\;\;0\;\;}$.
In this case, $x_p\in\la\ubrace{h_2}{m-1},y_{p-2m+1}\ra$.

(iii) Suppose $x_p\in Z^{p\bar0}_{2m+1}$ represents a class 
$[x_p]_{2m+1}\in E^{p\bar0}_{2m+1}$.
Let $x$ be given by (i) and let $\tilde x=x+x_{p+2m}$ with any choice of
$x_{p+2m}\in A^{p+2m}$.
Then $\bd_{\hhh^{(m+1)}}\tilde\xxx_p^{(m+1)}=\left(\substack{z_{p+2m+1}\\ 
\ooo_{p+1}^{(m)}}\right)$ for some $z_{p+2m+1}\in A^{p+2m+1}$, or equivalently,
${\hhh^{(m+1)}\;\;\tilde\xxx_p^{(m+1)}\choose\!\!\!\!\!\!0\quad\quad\;\;\;0}
\in\mat{m+2}_+(A)$ is a defining system (of $\ubrace{[h_2]}{m},[x_p]$).
Furthermore, $z_{p+2m+1}=\mu{\hhh^{(m+1)}\;\;\tilde\xxx_p^{(m+1)}\choose
\!\!\!\!\!\!0\quad\quad\;\;\;0}\in A^{p+2m+1}$ descends to a class in
$E^{p+2m+1,\bar0}_{2m+1}$ which is equal to $\dd_{2m+1}[x_p]_{2m+1}$.
\end{thm}

\noindent{\em Proof:}
(i) By the description of $Z^{p\bar q}_r$ in the Appendix,
$x_p\in Z^{p\bar0}_{2m+1}$ if and only if there is
$x=x_p+x_{p-2}+\cdots+x_{p+2m-2}\in A^p\op A^{p+2}\op\cdots\op A^{p+2m-2}$
such that $\dd x\in\Fi^{p+2m+1}A$.
The rest follows from Lemma~\ref{mat-bd}(i) and Corollary~\ref{mat-df}.\\
(ii) By the description of $B^{p\bar q}_r$ in the Appendix,
$x_p\in B^{p\bar0}_{2m+1}$ if and only if there is
$y=y_{p-2m+1}+y_{p-2m+3}+\cdots+y_{p-1}\in A^{p-2m+1}\op A^{p-2m+3}
\op\cdots\op A^{p-1}$ such that $\dd y-x_p\in\Fi^{p+3}A$.
The rest follows from Lemma~\ref{mat-bd}(i) and Corollary~\ref{mat-df}.\\
(iii) By the description of $\dd_r$ in the Appendix, $\dd_{2m+1}[x_p]_{2m+1}$
is represented by
$z_{p+2m+1}\in Z^{p+2m+1,\bar0}_{2m+1}\subset A^{p+2m+1}$ such that
$d\tilde x-z_{p+2m+1}\in\Fi^{p+2m+3}A$.
The rest follows from Lemma~\ref{mat-bd}(i) and Corollary~\ref{mat-df}.
\qed

Given $x_p$ in Theorem~\ref{spec-thm}(i), the element $x$ satisfying the
condition is not unique.
As in \S\ref{spec3}, this freedom does not exhaust all the defining systems
of $\ubrace{[h_2]}{m-1},[x_p]$.
Let $x_p+x'$, where $x'=x'_{p+2}+\cdots+x'_{p+2(m-1)}$, be another choice and
let $z'_{p+2m+1}=\mu{\hhh^{(m+1)}\;\;
\tilde\xxx'^{(m+1)}_p\choose\!\!\!\!\!\!0\quad\quad\;\;\;0}$,
where $\tilde x=x_p+x'+x'_{p+2m}$ for some $x'_{p+2m}\in A^{p+2m}$.
Then the difference
\begin{eqnarray*}
z'_{p+2m+1}-z_{p+2m+1}
\eq\mu{\hhh^{(m+1)}\;\;\tilde\xxx'^{(m+1)}_p-\tilde\xxx_p^{(m+1)}
  \choose\!\!\!\!\!\!0\quad\quad\quad\quad\quad\;\;\;0\quad\quad}        \\    
\eq\mu{\hhh^{(m)}\;\;\xxx'^{(m)}_{p+2}\choose\!\!0\quad\quad0\;\;}+
  \bb_1(x'_{p+2m}-x_{p+2m})
\end{eqnarray*}
descends to an element in the image of $\dd_{2m-1}$ which is a Massey product
of $\ubrace{[h_2]}{m-2},[x'_{p+2}-x_{p+2}]$. 
Therefore the class in $E_{2m}=E_{2m+1}$ remains the same.
This generalizes the discussion on the special case ($m=2$) on $\dd_5$ and
the triple Massey product in \S\ref{spec3}.

We also remark that the elements $x_{p+2},\dots,x_{p+2(m-1)}$ chosen in
Theorem~\ref{spec-thm} for a given $m\ge1$ can not be used recursively
without correction for higher values of $m$.
This phenomenon already appeared in the special case when $A$ is the de Rham
complex \cite{LLW}.
We now illustrate this fact for general $A_\infty$-algebras in a more concise
way.
Suppose the element $x_p$ in Theorem~\ref{spec-thm}(i) is actually in
$Z^{p\bar0}_{2m+3}$, then $\dd_{2m+1}[x_p]_{2m+1}$ given
by Theorem~\ref{spec-thm}(iii) is zero in $E^{p+2m+1,\bar0}_{2m+1}$, i.e., 
$z_{p+2m+1}\in B^{p+2m+1,\bar0}_{2m+1}$.
By Theorem~\ref{spec-thm}(ii), there exists 
$y'=y'_{p+2}+y'_{p+4}+\cdots+y'_{p+2m}$, where $y'_{p+2i}\in A^{p+2i}$
($1\le i\le m$) such that $\bd_{\hhh^{(m)}}\,\yyy'^{(m)}_{p+2}=
\left(\substack{z_{p+2m+1}\\ \ooo_{p+3}^{(m-1)}}\right)$.
Therefore
\[  \bd_{\hhh^{(m+1)}}\,\left(\tilde\xxx_p^{(m+1)}-
{\yyy'^{(m)}_{p+2}\choose0_p}\right)=\ooo_{p+1}^{(m+1)}.   \]
Thus when $m$ increases by $1$, the element $x$ in Theorem~\ref{spec-thm}(i)
should be replaced by
\[ x+x_{p+2m}-y'=x_p+(x_{p+2}-y'_{p+2})+\cdots+(x_{p+2m}-y'_{p+2m}).   \]

\subsection{Naturality of the spectral sequence}

If $\ff\colon A\to B$ is a morphism of $A_\infty$-algebras and $h\in A$ is
a twisting element, then so is $\FF_h(h)\in B$ and there is an induced
homomorphism $(\FF_h)_*\colon\Coh_h(A)\to\Coh_{\FF_h(h)}(B)$.
We want to describe explicitly the induced morphism of the spectral sequences.
If $h$ is of homogeneous degree $2$ as in the example in \S\ref{spec3} and
if $\ff$ is strict, then $\FF_h(h)=\ff_1(h)\in B$ is of homogeneous degree
$2$ as well.
In this case, the morphism of spectral sequences is induced by $\ff_1$.
It is easy to see, using the explicit formulas of $\dd_3$ and $\dd_5$ in
\S\ref{spec3}, that
\[  (\ff_1)_*\circ\dd^A_3=\dd^B_3\circ(\ff_1)_*,\quad
    (\ff_1)_*\circ\dd^A_5=\dd^B_5\circ(\ff_1)_*.   \]
In the general situation, we have

\begin{pro}\label{spec-f}
If $\ff\colon A\to B$ is a morphism of $A_\infty$-algebras and $h\in A$ is
a twisting element, then the cochain map $\FF_h\colon(A,\dd^A)\to(B,\dd^B)$,
where $\dd^A=\bd_h$ and $\dd^B=\bd_{\FF_h(h)}$, is a morphism of filtered
cochain complexes and hence there is a morphism of the spectral sequences
$({}^A\!E^{p\bar q}_r,\dd^A_r)\to({}^B\!E^{p\bar q}_r,\dd^B_r)$.
For all $m\ge0$, the homomorphisms $(\ff_1)_*\colon{}^A\!E^{p\bar0}_{2m+1}
\to{}^B\!E^{p\bar0}_{2m+1}$ are induced by $\ff_1\colon A\to B$ and 
$(\ff_1)_*\circ\dd^A_{2m+1}=\dd^B_{2m+1}\circ(\ff_1)_*$.
In particular, we have
$(\ff_1)_*\colon{}^A\!E^{p\bar0}_\infty\to{}^B\!E^{p\bar0}_\infty$.
\end{pro}

\noindent{\em Proof:}
The cochain map $\FF_h$ clearly preserve the filtrations.
On the graded components $\Gr^pA=A^p$, it is simply $\ff_1\colon A^p\to B^p$.
The rest is self-evident by the discussion of naturality in the Appendix.
\qed

The result that the morphism of the spectral sequence is induced by $\ff_1$
alone is consistent with the naturality of the higher Massey products in
\S\ref{nat-m}.
We note that ${}^A\!E^{p\bar0}_\infty$ are the graded components of
$\Coh^{\bar p}_h(A)$.
Although the homomorphisms on $E_\infty$ depend solely on $\ff_1$, the induced
map $(\FF_h)_*\colon\Coh_h(A)\to\Coh_{\FF_h(h)}(B)$ does depend on all
$\ff_n$ for $n\ge1$.

We apply Proposition~\ref{spec-f} to a number of cases.

First, suppose $\qq\colon\Coh(A)\to A$ is the quasi-isomorphism in
\S\ref{sec-ha}.
Denote by $(\bar E^{p\bar q}_r,\bar\dd_r)$ the spectral sequence associated
to the $\Z$-graded $A_\infty$-algebra $\Coh(A)$.
Since $\bar\bb_1=0$, we have $\bar E^{p\bar0}_3=\bar E^{p\bar0}_2=
\bar E^{p\bar0}_1=\bar E^{p\bar0}_0=\Coh^p(A)$.
The map $\qq_1\colon\bar E^{p\bar0}_0=\Coh^p(A)\to E^{p\bar0}_0=A^p$ is a
quasi-isomorphism of cochain complexes and hence the induced homomorphisms
$(\qq_1)_*\colon\bar E^{p\bar0}_r\to E^{p\bar0}_r$ are the identity maps 
for all $r\ge2$.
With the relation between the higher differentials and the higher Massey
products (Theorem~\ref{spec-thm}(iii)), this is consistent with the relation
of (higher) Massey products with the $A_\infty$-structure on $\Coh(A)$
(Proposition~\ref{ha-mas}).

Next, if as in \S\ref{hh'}, $h,h'\in A$ are two twisting elements
that are homotopic through $c$, then by Proposition~\ref{trans}(i)
and Corollary~\ref{isom-c}(i), there is a cochain map
$\uppsi_c\colon(A,\bd_h)\to(A,\bd_{h'})$,
which preserves the filtration on $A$.
Moreover, $\uppsi_c$ is the identity map on the graded components.
So the induced isomorphism identifies the two spectral sequences
$(E^{p\bar q}_r,\dd_r)=(E'^{p\bar q}_r,\dd'_r)$, which converge
to $\Coh_h(A)$ and $\Coh_{h'}(A)$, respectively.
We note that although $\Coh_h(A)$ and $\Coh_{h'}(A)$ have identical graded
components, the isomorphism between the total spaces is induced by $\uppsi_c$,
not the identity map.

Finally, if $\ff,\gg\colon A\to B$ are two morphisms of $A_\infty$-algebras
that are homotopic through $\hh$, then $\ff_1,\gg_1\colon A\to B$ are
cochain maps that are homotopic through $\hh_1$. 
Consequently, the morphisms on the spectral sequence ${}^A\!E^{p\bar q}_r$ 
induced by $\ff$ and $\gg$ are identical when $r\ge2$.
This is consistent with Theorem~\ref{fgh} and Corollary~\ref{fgh-isom} since
$\uppsi_{\HH_h(h)}$ induces the identity isomorphism on spectral sequences.

\section*{Appendix.\ Spectral sequences}

We summarize some aspects of spectral sequences used in Section~5.
We refer the readers to numerous discussions in the literature (for example
\cite{CE,Mc}) although the version we present here is somewhat different.
We will treat the cohomological spectral sequences only as the homological ones
are parallel.

Let $(C^\bullet,\dd)$ be a ($\Z$-graded) cochain complex with a filtration
$\Fi$, that is, we have a descending $\dd$-invariant sequence 
\[ \cdots\supset\Fi^pC^\bullet\supset\Fi^{p+1}C^\bullet\supset\cdots \]
and thus each $(\Fi^pC^\bullet,\dd)$ is a cochain complex.
Let $\Gr^pC^\bullet=\Fi^pC^\bullet/\Fi^{p+1}C^\bullet$ be the graded
components; each $\Gr^pC^\bullet=\Fi^pC^\bullet/\Fi^{p+1}C^\bullet$ is also
a cochain complex whose coboundary operator will also be denoted by $\dd$.
Let
\[ \tilde Z\pqr=\Fi^pC^{p+q}\cap\dd^{-1}(\Fi^{p+r}C^{p+q+1}),\quad
   \tilde B\pqr=\Fi^pC^{p+q}\cap\dd(\Fi^{p-r+1}C^{p+q-1})  \]
(which are usually denoted by $Z\pqr,B\pqr$, respectively) and
\[ E\pqr=\tilde Z\pqr/(\tilde B\pqr+\tilde Z\opqr). \]
Then the coboundary operator $\dd$ induces the higher differentials
$\dd_r\colon E\pqr\to E\pqrr$ with $\dd_r^2=0$, and $E_{r+1}^{pq}$ is the
cohomology of $(E\pqr,\dd_r)$.
The filtration is exhaustive and weakly convergent, then the spectral sequence
converges to the cohomology groups (cf.\ Theorem~3.2 of \cite{Mc}).

We reserve the notations $Z\pqr,B\pqr$ for the new spaces
\[ Z\pqr=\tilde Z\pqr/\tilde Z\opqr,\quad B\pqr=
         \tilde B\pqr/(\tilde B\pqr\cap\tilde Z\opqr). \]
Then $E\pqr=Z\pqr/B\pqr$.
There are natural isomorphisms (through which we identify)
\[ Z\pqr=\set{x\in\Gr^pC^{p+q}}{\exists\tilde x\in\Fi^pC^{p+q}
              \mbox{ such that }x=\tilde x+\Fi^{p+1}C^{p+q},\;
              \dd\tilde x\in\Fi^{p+r}C^{p+q+1}}                       \]
and
\[ B\pqr=\set{x\in\Gr^pC^{p+q}}{\exists\tilde y\in\Fi^{p-r+1}C^{p+q-1}
              \mbox{ such that }\dd\tilde y\in\Fi^pC^{p+q},\;
              x=\dd\tilde y+\Fi^{p+1}C^{p+q}}.                      \]
For example, $Z_0^{pq}=\Gr^pC^{p+q}$, $B_0^{pq}=0$ and
$E_0^{pq}=\Gr^pC^{p+q}$ whereas 
\[ Z_1^{pq}=\ker(\dd\colon\Gr^pC^{p+q}\to\Gr^pC^{p+q+1}),\quad
   B_1^{pq}={\mathrm{im}}(\dd\colon\Gr^pC^{p+q-1}\to\Gr^pC^{p+q}) \]
and $E_1^{pq}=H^{p+q}(\Gr^pC^\bullet,\dd)$ as expected.
We have a sequence of inclusions
\[ 0=B_0^{pq}\subset\cdots\subset B\pqr\subset B_{r+1}^{pq}
  \subset\cdots\subset Z_{r+1}^{pq}\subset Z\pqr\subset\cdots
  \subset Z_0^{pq}=\Gr^pC^{p+q}.                      \]
Finally, if the filtration is weakly convergent, then
\[ Z_\infty^{pq}=\set{x\in\Gr^pC^{p+q}}{\exists\tilde x\in\Fi^pC^{p+q}
      \mbox{ such that }x=\tilde x+\Fi^{p+1}C^{p+q},\;\dd\tilde x=0},    \]
\[ B_\infty^{pq}=\set{x\in\Gr^pC^{p+q}}{\exists\tilde y\in C^{p+q-1}
              \mbox{ such that }\dd\tilde y\in\Fi^pC^{p+q},\;
              x=\dd\tilde y+\Fi^{p+1}C^{p+q}},                      \]
and $E_\infty^{pq}=Z_\infty^{pq}/B_\infty^{pq}$ are the
graded components of the cohomology groups $\Coh(C^\bullet,\dd)$.

We now describe the higher differentials $\dd_r\colon E\pqr\to E\pqrr$.
If a class $[x]\in E\pqr=Z\pqr/B\pqr$ is represented by
$x\in Z\pqr\subset\Gr^pC^{p+q}$, choose a lifting
$\tilde x\in\Fi^pC^{p+q}$ of $x$.
Then $\dd\tilde x\in\Fi^{p+r}C^{p+q+1}$ and we denote by 
$z\in\Gr^{p+r}C^{p+q+1}$ its graded component.
In fact, $z\in Z\pqrr$ (since if we choose $\tilde z=\dd\tilde x$, then
$\dd\tilde z=0\in\Fi^{p+2r}C^{p+q+2}$) and $\dd_r$ is given by
\[  \dd_r[x]=[z]\in E\pqrr.  \]
We check that the result is independent of the choices made.
If $\tilde x'\in\Gr^pC^{p+q}$ is another lifting of $x$, then
$\tilde x'-\tilde x\in\Fi^{p+1}C^{p+q}$ and 
$\dd(\tilde x'-\tilde x)\in\Fi^{p+r}C^{p+q+1}$.
Since the graded component $z'\in\Gr^{p+r}C^{p+q}$ of $\dd\tilde x'$ differs
from $z$ by an element in $B\pqrr$, they descend to the same class in $E\pqrr$.
On the other hand, if $[x]=0$ or $x\in B\pqr$, then we can choose
$\tilde x=\dd\tilde y$, where $\tilde y\in\Fi^{p-r+1}C^{p+q-1}$.
Thus $z=\dd\tilde x=0$ and $\dd_r[x]=0$, which is required by consistency.

Let $(C'^\bullet,\dd')$ be another ($\Z$-graded) cochain complex with a
filtration, also denoted by $\Fi$.
Suppose $\ff\colon(C^\bullet,\dd)\to(C'^\bullet,\dd')$ is a morphism of
filtered cochains, i.e., $\ff\circ\dd=\dd'\circ\ff$ and
$\ff\colon\Fi^pC^\bullet\to\Fi^pC'^\bullet$.
Then there is an induced morphism of the associated spectral sequences
$(E\pqr,\dd_r)\to(E'^{pq}_r,\dd'_r)$ (see, for example, Theorem~3.5 of
\cite{Mc}).
More concretely, $\ff$ induces cochain maps (using the same notation) 
$\ff\colon\Gr^pC^\bullet\to\Gr^pC'^\bullet$.
It can be shown that
$\ff(Z\pqr)\subset Z'^{pq}_r$, $\ff(B\pqr)\subset B'^{pq}_r$
and hences $\ff$ induces a map $\ff_*\colon E\pqr\to E'^{pq}_r$.
Furthermore, we have $\ff_*\circ\dd_r=\dd'_r\circ\ff_*$, which can also
be seen from the above explicit description of $\dd_r$.
This gives a morphism between the two spectral sequences.

\medskip

\noindent{\bf Acknowledgments.} 
S.W.\ thanks Oklahoma State University and W.L.\ thanks University of
Colorado and CUHK/IMS, where part of the work was done, for hospitality.
We thank J.\ Palmieri and A.\ Prout\'e for discussions and J.\ Stasheff
for comments on the manuscript.

\end{document}